\documentclass[a4paper,12pt,reqno]{amsart}
\usepackage{graphics,hyperref}
\newtheorem{lem}{Lemma}

\newtheorem{theo}{Theorem}

\newtheorem{cor}{Corollary}
\newtheorem{defi}{Definition}
\numberwithin{equation}{section}

\newcommand{\sgn}{\operatorname{sgn}}
\newcommand{\norm}{\operatorname{norm}}
\newcommand{\SP}{\operatorname{SP}}
\newcommand{\qbinom}[2]{\left[ #1 \atop #2 \right]}

\begin{document}

\title[Bender-Knuth (ex-)Conjecture]{Another refinement of the Bender-Knuth (ex-)Conjecture}

\author[Ilse Fischer]{\box\Adr}

\newbox\Adr
\setbox\Adr\vbox{ \centerline{ \large Ilse Fischer} \vspace{0.3cm}
\centerline{Institut f\"ur Mathematik, Universit\"at Klagenfurt,}
\centerline{Universit\"atsstrasse 65-67, A-9020 Klagenfurt, Austria.}
\centerline{E-mail: {\tt Ilse.Fischer@uni-klu.ac.at}}}

\maketitle

\begin{abstract}
We compute the generating function of column-strict plane
partitions with  parts in $\{1,2,\ldots,n\}$, at most $c$
columns, $p$ rows of odd length and $k$ parts equal
to $n$. This refines both, Krattenthaler's~\cite{kratt} and the
author's~\cite{fischer} refinement of the Bender-Knuth
(ex-)Conjecture. The result is proved by an extension of the
method for proving polynomial enumeration formulas which was
introduced by the author in~\cite{fischer} to $q$-quasi-polynomials.
\end{abstract}

\section{Introduction}

Let $\lambda=(\lambda_{1},\lambda_{2},\dots,\lambda_{r})$ be a
partition, i.e. $\lambda_i \in \mathbb{Z}$ and $\lambda_1 \ge
\lambda_2 \ge \ldots \ge \lambda_r \ge 0$. A {\it strict plane
partition} of shape $\lambda$ is an array ${\Pi} = (\pi_{i,j})_{1 \le i \le r,
1 \le j \le \lambda_{i}}$ of non-negative integers such that the rows
are weakly decreasing and the columns are strictly decreasing. For
instance
$$
\begin{array}{ccccc}
7 & 6 & 5 & 5 & 2 \\
5 & 4 & 2 & 2 &   \\
4 & 2 &   &   &   \\
2 & 1 &   &   &
\end{array}
$$
is a strict plane partition of shape $(5,4,2,2)$. The {\it norm}
$n({\Pi})$ of a strict plane partition is defined as the sum of its
parts and ${\Pi}$ is said to be a strict plane partition of the
non-negative integer $n({\Pi})$. Thus $47$ is the norm of our example.
Strict plane partitions and closely related objects
have been enumerated subject to a variety of different
constraints. In \cite[p.50]{benderknuth} Bender and Knuth had
conjectured that the generating function with respect to the norm
of strict plane partitions with at most $c$ columns and parts in
$\{1,2,\dots,n\}$ is equal to
$$
\sum q^{n(\pi)} = \prod_{i=1}^n \frac{[c+i;q]_i}{[i;q]_i},
$$
where $[n;q] = 1 + q + \dots + q^{n-1}=(1-q^{n})/(1-q)$ and $[a;q]_n =
\prod_{i=0}^ {n-1} [a+i;q]$. This conjecture was proved by
Andrews~\cite{andrews}, Gordon~\cite{gordon}, Macdonald~\cite[Ex.
19, p.53]{macdonald} and Proctor~\cite[Prop. 7.2]{proctor}. For
related papers, which mostly include generalizations of the
Bender-Knuth (ex-)Conjecture, see~\cite{des1, des2, fischer,
kadell, kratt, proctor2, stem2}.

In particular, Krattenthaler~\cite{kratt} computed the generating function of
strict plane partitions with parts in $\{1,2,\ldots,n\}$, at
most $c$ columns and $p$ rows of odd length. On the
other hand the author~\cite{fischer} computed the generating function of strict plane partitions
with parts in $\{1,2,\ldots,n\}$, at most $c$ columns and $k$
parts equal to $n$. In this paper we refine
these two results. The main result is the following.

\begin{theo}
\label{main} The generating function of strict plane partitions
with parts in $\{1,2,\ldots,n\}$, at most $c$ columns, $p$ rows of
odd length and $k$ parts equal to $n$ is given by
\begin{multline*}
 M_{n,c,p} [k+1;q]_{n-1} [k-c-n+1;q]_{n-1} q^k +
L_{n,c,p} \Bigg ((-1)^{k} q^{n k} + (-1)^n q^{(n-1)(2c+n)/2+k} \\ \times
\sum_{i=1}^{n-1}
(-1)^{c} q^{\binom{i}{2}} \frac{[k+1;q]_{n-1}
  [k-c-i+1;q]_{i-1} [k-c-n+1;q]_{n-i-1}}{[1;q]_{i-1} [1;q]_{n-1-i}
  [c+i+1;q]_{n-1}} \\
- q^{\binom{i}{2}} \frac{[k+1;q]_{i-1}
  [k+i+1;q]_{n-i-1} [k-c-n+1;q]_{n-1}}{ [1;q]_{i-1}
  [1;q]_{n-1-i} [c+i+1;q]_{n-1}}  \Bigg)
\end{multline*}
where
{\small \begin{equation*}
L_{n,c,p} =
\begin{cases}
\left(
\frac{q^{\binom{p+1}{2}} {\qbinom{n-1}{p}} [c;q]}{[c+p;q]_n} -
\frac{q^{\binom{p}{2}} {\qbinom{n-1}{p-1}} [c+2n;q]}{[c+p+1;q]_n}
\right) \frac{[c+1;q]_{n-1} [1;q]_{n-1}}{2}
\prod\limits_{i=1}^{n-1} \frac{[c+2 i + 1;q]_{n-i}}{[2 i;q]_{n-i} [2 i;q]} & 2|c \\
\left( q^{\binom{p+1}{2}}
{\qbinom{n-1}{p}} - q^{\binom{p}{2}} {\qbinom{n-1}{p-1}} \right) \frac{[1;q]_{n-1}}{2}
\prod\limits_{i=1}^{n-1} \frac{[c+2 i;q]_{n-i}}{[2 i;q]_{n-i} [2 i;q]}
& 2 \not| c
\end{cases}
\end{equation*}}
and
{\small
\begin{multline*}
M_{n,c,p}= \frac{(-1)^{n-1} q^{(n-1)(2c+n)/2}}{[1;q]_{n-2}} \\
\times
\begin{cases}
\left(
\frac{q^{\binom{p+1}{2}} \qbinom{n-1}{p} [c;q]}{[c+p;q]_{n}}
\left( \frac{1}{[n-1;q]} -
\frac{[c+2n-1;q]}{[c+n;q] [2n-2;q]} \right)  +
\frac{q^{\binom{p}{2}} {\qbinom{n-1}{p-1}}}{[c+p+1;q]_n}
\frac{[c+2n-1;q]_2}{[c+n;q] [2n-2;q]} \right) \prod\limits_{i=1}^{n-1}
\frac{[c+2i;q]_{n-i}}{[2i;q]_{n-i}} & 2|c \\
\left( q^{\binom{p+1}{2}+n-1} {\qbinom{n-1}{p}}
 + q^{\binom{p}{2}} {\qbinom{n-1}{p-1}} \right) \frac{1}{[2n-2;q]} \prod\limits_{i=1}^{n-1}
\frac{[c+2i+1;q]_{n-i-1}}{[2i;q]_{n-i}} & 2\not|c
\end{cases}.
\end{multline*}}
\end{theo}

\medskip

In these formulas the notion of the $q$-binomial coefficient is used. It is
defined as follows.
\begin{equation*}
\qbinom{n}{k} =
\begin{cases}
\frac{[n-k+1;q]_k}{[1;q]_k}  & \text{if $0 \le k \le n$} \\
0   & \text{otherwise}
\end{cases}
\end{equation*}
At the end of Section~\ref{compute} we show that Theorem~\ref{main} implies
Krattenthaler's and the author's refinement of the Bender-Knuth (ex-)Conjecture
and with this the Bender-Knuth (ex-)Conjecture itself.

Our method for proving Theorem~\ref{main} is an extension of the
method for proving polynomial enumeration formulas we have
introduced in \cite{fischer}. It is interesting to note that this
elementary method avoids the use of determinants completely, which
is quite unusual in the field of plane partition enumeration. The
method is divided into the following three steps.

\begin{enumerate}
\item {\bf Extension of the combinatorial interpretation.} It only makes sense
to ask for the number of strict plane partitions with parts in
$\{1,2,\ldots,n\}$, at most $c$ columns and $k$ parts equal to $n$
if $k \in \{0,1,\ldots,c\}$. This is because all $n$'s must be in
the first row of the strict plane partition by the columnstrictness. In the first step we
find a combinatorial extension of these strict plane partitions to
arbitrary integers $k$, i.e. we find new objects indexed by an
arbitrary integer $k$ which are in bijection with strict plane
partitions with $k$ parts equal to $n$ if $k \in
\{0,1,\ldots,c\}$.

\item {\bf The extending objects are enumerated by a $q$-quasi-polynomial in
$k$.} With the help of a simple recursion we show that the
extending objects are enumerated by a $q$-quasi-polynomial (see
Definition~\ref{q-quasi-polynomial}). Moreover the degree of the
$q$-quasi-polynomial is computed.

\item {\bf Exploring properties of the $q$-quasi-polynomial that determine it
uni\-quely.} A ($q$-quasi-)polynomial is determined by a finite
number of properties such as zeros or other evaluations. In the
last step we derive enough properties of the $q$-quasi-polynomial
in order to compute it using the degree estimation from the
previous step.
\end{enumerate}

Note that this article contains two types of extensions of the
method for proving polynomial enumeration formulas presented in
\cite{fischer}. Firstly, the method is extended to
$q$-quasi-polynomials, see Definition~\ref{q-quasi-polynomial}.
More remarkable is, however, the following extension: In
\cite{fischer} we have described a method that is applicable to
polynomial enumeration formulas that factorize into distinct
linear factors over $\mathbb{Z}$. There the
``properties'' in the third step are just the integer zeros together
with one  (easy to compute) non-zero evaluation. (Thus the third
step was entitled ``Exploring natural linear factors''.) In this
article we demonstrate that the lack of enough integer zeros can
be compensated by other properties of the ($q$-quasi-)polynomial.

\smallskip

The paper is organized as follows. In Section~\ref{ext} we give
the combinatorial extension of strict plane partitions as proposed
in Step 1. In Section~\ref{qquasi} we introduce the notion of
$q$-quasi-polynomials and establish the properties needed in this
paper. In Section~\ref{isqquasi} we show that the
generating function of strict plane partitions which is under consideration
in this paper is a $q$-quasi-polynomial and we compute its degree (Step 2). In
Section~\ref{properties} we deduce enough properties of the
$q$-quasi-polynomial in order to compute it (Step 3). In
Section~\ref{compute} we perform the (complicated) computation and
in Section~\ref{q-sums} we derive some $q$-summation formulas
which are needed in the computation.

\medskip

Throughout the whole article we use the extended definition of
the summation symbol, namely,
\begin{equation}
\label{ext-sum}
\sum_{i=a}^{b} f(i) =
\begin{cases}
f(a) + f(a+1) + \dots + f(b) & \text{if $a \le b$} \\
0 & \text{if $b=a-1$} \\
- f(b+1) - f(b+2) - \dots - f(a-1) & \text{if $b+1 \le a-1 $}
\end{cases}.
\end{equation}
This assures that for any polynomial $p(X)$ over an arbitrary
integral domain $I$ containing $\mathbb{Q}$ there exists a unique
polynomial $q(X)$ over $I$ such that $\sum_{x=0}^y p(x) = q(y)$
for all integers $y$. We usually write $\sum_{x=0}^ y p(x)$ for
$q(y)$.

\section{Extension of the combinatorial interpretation}
\label{ext}

In this section we establish the combinatorial extension of strict
plane partitions with parts in $\{1,2,\ldots,n\}$, at most $c$
columns and $k$ parts equal to $n$ to arbitrary integers $k$. This
extension was already introduced in \cite[Section~2]{fischer}. We
repeat it here in less detail.

Let $r,n,c$ be integers with $0 \le r \le n$. A {\it generalized
$(r,n,c)$ Gelfand-Tsetlin-pattern} (for short: $(r,n,c)$-pattern)
is an array $(a_{i,j})_{1 \le i \le r+1, i-1 \le j \le n+1}$ of
integers with

\begin{enumerate}

\item $a_{i,i-1}=0$ and $a_{i,n+1}=c$,

\item if $a_{i,j} \le a_{i,j+1}$ then $a_{i,j} \le a_{i-1,j} \le a_{i,j+1}$

\item if $a_{i,j} > a_{i,j+1}$ then $a_{i,j} > a_{i-1,j} > a_{i,j+1}.$

\end{enumerate}

The norm of an $(r,n,c)$-pattern is defined as the sum of its
parts, where the first and the last part of each row is omitted. A
$(3,6,c)$-pattern for example is of the form

\vspace{5mm}

\begin{tabular}{ccccccccccccccc}
 &   & & $0$ & & $a_{4,4}$ & & $a_{4,5}$ & & $a_{4,6}$ & & $c$ &  &  & \\
 & &$0$ &   & $a_{3,3}$ & & $a_{3,4}$ & & $a_{3,5}$ & & $a_{3,6}$ & & $c$ & & \\
 & $0$ & & $a_{2,2}$ & & $a_{2,3}$ & & $a_{2,4}$ & & $a_{2,5}$ & & $a_{2,6}$ & & $c$ & \\
$0$ & & $a_{1,1}$ & & $a_{1,2}$ & & $a_{1,3}$ & & $a_{1,4}$ & &
$a_{1,5}$ & & $a_{1,6}$ & & $c$,
\end{tabular}

\vspace{5mm}

\noindent such that every entry not in the top row is between its
northwest neighbour $w$ and its northeast neighbour $e$, if $w \le
e$ then weakly between, otherwise strictly between. Thus

\vspace{5mm}

\begin{center}
\begin{tabular}{ccccccccccccccc}
 &   & & $0$ & & $3$ & & $-5$ & & $10$ & & $4$ &  &  & \\
 & &$0$ &   & $2$ & & $-2$ & & $3$ & & $8$ & & $4$ & & \\
 & $0$ & & $2$ & & $-1$ & & $2$ & & $4$ & & $7$ & & $4$ & \\
$0$ & & $0$ & & $0$ & & $1$ & & $2$ & & $5$ & & $6$ & & $4$
\end{tabular}
\end{center}

\vspace{5mm}

\noindent is an example of an $(3,6,4)$-pattern. Note that a
generalized $(n-1,n,c)$ Gelfand-Tsetlin-pattern $(a_{i,j})$ with
$0 \le a_{n,n} \le c$ is what is said to be a
Gelfand-Tsetlin-pattern with $n$ rows and parts in
$\{0,1,\ldots,c\}$, see~\cite[p. 313]{stanley} or
\cite[(3)]{gelfand} for the original reference.
(Observe that $0 \le a_{n,n} \le c$ implies that the third condition in the
definition of a generalized Gelfand-Tsetlin-pattern never applies.)
The following
correspondence between Gelfand-Tsetlin-patterns and strict plane
partitions is crucial for our paper.

\begin{lem}
\label{bijection} There exists a norm-preserving bijection between
Gelfand-\-Tsetlin-\-patterns  with $n$ rows, parts in
$\{0,1,\ldots,c\}$ and fixed $a_{n,n}=k$, and strict plane
partitions with parts in $\{1,2,\dots,n\}$, at most $c$ columns
and $k$ parts equal to $n$. In this bijection
$(a_{1,n},a_{1,n-1},\ldots,a_{1,1})$ is the shape of the strict
plane partition.
\end{lem}

{\it Proof.} Given such a Gelfand-\-Tsetlin-\-pattern, the
corresponding strict plane partition is such that the shape filled
by parts greater than $i$ corresponds to the partition given by
the $(n-i)$-th row (the top row being the first row) of the
Gelfand-Tsetlin-pattern, where the first and the last part of the
row in the pattern are omitted. Thus the strict plane partition in
the introduction corresponds to the following Gelfand-Tsetlin
pattern (first and last parts in the rows are omitted).
\begin{center}
\begin{tabular}{ccccccccccccccccc}
  &   &   &   &   &   &   &   & 1 &   &   &   &   &   &   &   & \\
  &   &   &   &   &   &   & 0 &   & 2 &   &   &   &   &   &   & \\
  &   &   &   &   &   & 0 &   & 1 &   & 4 &   &   &   &   &   & \\
  &   &   &   &   & 0 &   & 1 &   & 2 &   & 4 &   &   &   &   & \\
  &   &   &   & 0 &   & 0 &   & 1 &   & 2 &   & 4 &   &   &   & \\
  &   &   & 0 &   & 0 &   & 1 &   & 2 &   & 4 &   & 5 &   &   & \\
  &   & 0 &   & 0 &   & 0 &   & 1 &   & 2 &   & 4 &   & 5 &   &
\end{tabular}
\end{center}
\qed

\medskip

Therefore it suffices to compute the generating function with
respect to the norm of $(n-1,n,c)$-patterns with fixed
$a_{n,n}=k$, $0 \le k \le c$, and where exactly $p$ values of
$a_{1,1},a_{1,2},\ldots,a_{1,n}$ are odd. However,
$(n-1,n,c)$-patterns are defined for all $a_{n,n} \in \mathbb{Z}$
and thus we have established the combinatorial extension apart
from the following technical detail. That is that we actually have
to work with a signed enumeration if $a_{n,n} \notin
\{0,1,\ldots,c\}$. Therefore we define the sign of a pattern.

A pair  $(a_{i,j}, a_{i,j+1})$ with $a_{i,j} > a_{i,j+1}$ and $i
\not= 1$ is called an {\it inversion} of the $(r,n,c)$-pattern and
$(-1)^{\# \, \text{of inversions}}$ is said to be the {\it sign}
of the pattern, denoted by $\sgn(a)$. The $(3,6,4)$-pattern in the
example above has altogether $6$ inversions and thus its sign is
$1$. We define the following generating function
$$
F_q(r,n,c,p;k_{1},k_{2},\dots,k_{n-r}) = \left( \sum_{a} \sgn(a)
q^{\norm(a)} \right) / q^{k_1+k_2+\ldots+k_{n-r}},
$$
where the sum is over all $(r,n,c)$-patterns $(a_{i,j})$ with top
row defined by $k_{i} = a_{r+1,r+i}$ for $i=1,\dots,n-r$ and such
that exactly $p$ of $a_{1,1},a_{1,2},\ldots,a_{1,n}$ are odd. It
is crucial that for $0 \le k \le c$ $F_{q}(n-1,n,c,p;k) \, q^{k}$ is the generating function of
$(n-1,n,c)$-patterns with $a_{n,n}=k$ and where
exactly $p$ of $a_{1,1},a_{1,2},\ldots,a_{1,n}$ are odd.
This is because an $(n-1,n,c)$-pattern with $0 \le a_{n,n} \le c$ has no
inversions. Thus $F_{q}(n-1,n,c,p;k)$ is the quantity we want to
compute. It has the advantage that it is well defined for all
integers $k$, whereas our original enumeration problem was only
defined for $0 \le k \le c$.

\section{$q$-quasi-polynomials and their properties}
\label{qquasi}

In the following let $R$ be a ring containing
$\mathbb{C}$. A {\it quasi-polynomial} (see~\cite[page
210]{stanley1}) in the variables $X_1, X_2, \ldots, X_n$ over $R$ is
an expression of the form
$$
\sum_{(m_1,m_2, \ldots, m_n) \in \mathbb{Z}^n, m_i \ge 0}
c_{m_1,m_2,\ldots,m_n} (X_1,X_2,\ldots,X_n) \, X_1^{m_1} X_2^{m_2}
\cdots X_n^{m_n},
$$
where $(X_1,X_2,\ldots,X_n) \to
c_{m_1,m_2,\ldots,m_n}(X_1,X_2,\ldots,X_n)$ are periodic functions
on $\mathbb{Z}^n$ taking values in $R$, that is there exists an integer $t$ with
$$
c_{m_1,m_2,\ldots,m_n} (k_1,\ldots,k_i,\ldots,k_n) =
c_{m_1,m_2,\ldots,m_n} (k_1,\ldots,k_{i}+t,\ldots,k_n)
$$
for all $(k_{1},\ldots,k_{n}) \in \mathbb{Z}^{n}$ and $i$, and
almost all $c_{m_1,\ldots,m_n}(X_1,\ldots,X_n)$ are zero. Let
$(m_1,\ldots,m_n)$ be with $c_{m_1,\ldots,m_n} (X_1,\ldots,X_n)
\not= 0$ such that $m_1  + \ldots + m_n$ is maximal. Then $m_1 +
\ldots + m_n$ is said to be the degree of the quasi-polynomial.
(The zero-quasi-polynomial is said to be of degree $- \infty$.)
The smallest common period of all
$c_{m_1,\ldots,m_n}(X_1,\ldots,X_n)$ is said to be the period of
the quasi-polynomial.  (In this paper we only deal with
$q$-quasi-polynomials of period 1 or 2.) In \cite[Section 6]
{fischer} we have defined $q$-polynomials. The following
definition of $q$-quasi-polynomials is the merge of these two
definitions. In this definition let $R_q$ denote the ring of
quotients  with elements from $R[q]$ in the numerator and elements
from $\mathbb{C}[q]$ in the denominator.

\begin{defi}
\label{q-quasi-polynomial} A $q$-quasi-polynomial over $R$ in
$X_1,X_2,\ldots,X_n$ is a quasi-polynomial over $R_q$ in $q^{X_1},
q^{X_2},\ldots, q^{X_n}$. Let $R_{qq}[X_1,X_2,\ldots,X_n]$ denote
the ring of these $q$-quasi-polynomials.
\end{defi}

Observe that $R_{qq}[X_1,\ldots,X_n]$ is the ring of
$q$-quasi-polynomials in $X_i$ over
$$R_{qq}[X_1,\ldots,X_{i-1},X_{i+1},\ldots,X_n].$$ (Thus it would have
been possible to define $R_{qq}[X_1,\ldots,X_n]$ inductively
with respect to $n$.) We define $[X;q]=(1-q^{X})/(1-q)$ and $[X;q]_{n}=\prod_{i=0}^{n-1}
[X+i;q]$. Observe that $$[X_{1};q]_{m_{1}} [X_{2};q]_{m_{2}}
\cdots [X_{n};q]_{m_{n}},
$$
with $(m_{1},m_{2},\ldots,m_{n}) \in \mathbb{Z}^{n}$ and $m_{i}
\ge 0$, is a basis of the $q$-quasi-polynomials over the periodic
functions.

The following two properties of polynomials were crucial for our
method for proving polynomial enumeration formulas which we have
introduced in \cite{fischer}. Since we want to extend our method to
$q$-quasi-polynomials, we have to find $q$-quasi-analogs of these
properties.

\begin{enumerate}
\item If $p(X)$ is a polynomial over $R$, then there exists a (unique)
polynomial $r(X)$ with $\deg r = \deg p +1$ and
$$
\sum_{x=0}^y p(x) = r(y)
$$
for every integer $y$.

\item If $p(X)$ is a polynomial over $R$ and $a$ is a zero of $p(X)$, then there exists a polynomial
$r(X)$ over $R$ with
$$
p(X) = (X-a) r(X).
$$
\end{enumerate}

Regarding the first property we show the following for $q$-quasi-polynomials.

\begin{lem}
\label{q-quasi-sum}
 Let $p(X)$ be a $q$-quasi-polynomial in $X$ over $R$ with degree $d$ and period $t$.
Then $\sum_{x=1}^y p(x) q^x$ is a $q$-quasi-polynomial over $R$ in $y$ with
degree at most $d+1$ and  period at most $t$.
\end{lem}

In order to prove this lemma we need a definition and another
lemma.

\begin{defi} Let $\rho \to f(\rho)$ be a function. Then the $q$-differential-operator $\frac{d}{d_q \rho}$
is defined as follows
$$
\frac{d}{d_q \rho} f(\rho) = \frac{f(q \, \rho)- f(\rho)}{\rho
(q-1)}.
$$
With $\frac{d}{d_q \rho^n}$ we denote the $n$-fold application of
the operator.
\end{defi}

Observe that for a laurent polynomial we have
\begin{equation}
\label{q-diff-poly}
\frac{d}{d_q \rho} \sum_{i=b}^c a_i \rho^i =
\sum_{i=b}^c [i;q] a_i \rho^{i-1}.
\end{equation}
Note that this is also true if $b > c$.

\begin{lem}
\label{identity}
$$\sum_{x=0}^y [x;q]_n q^{x} \sigma^{x-1} = \frac{d}{d_q \sigma^n} \left( \frac{\sigma^{n-1} ((\sigma q)^{y+1} -
1)}{(\sigma q-1)} \right)$$
\end{lem}

{\it Proof of Lemma~\ref{identity}.} By \eqref{q-diff-poly} we have the
following identity.
$$
\sum_{x=0}^y [x;q]_n q^{x} \sigma^{x-1} = \frac{d}{d_q \sigma^n}
\left( \sum_{x=0}^y q^x \sigma^{x+n-1} \right).
$$
The assertion now follows from
$$
\sum_{x=0}^y  q^x \sigma^{x+n-1} = \frac{\sigma^{n-1} ((\sigma q
)^{y+1} - 1)}{(\sigma q-1)}.
$$
\qed

{\it Proof of Lemma~\ref{q-quasi-sum}.} Suppose $p(X)$ is a
$q$-quasi-polynomial with period $t$. Let $\rho \in \mathbb{C}$ be a primitive
$t$-th root of unity. Then $p(X)$ can be expressed as follows
$$
p(X) = p_0(X) + \rho^X p_1(X) + \rho^{2 X} p_2(X) + \ldots +
\rho^{(t-1) X} p_{t-1}(X),
$$
where $p_i(X)$ are $q$-polynomials, i.e. $q$-quasi-polynomials
with period $1$. Suppose $d$ is the degree of $p(X)$. Then, for
every $i$, we have
$$
p_i(X) = \sum_{j=0}^{d} a_{i,j} \, [X;q]_j,
$$
where $a_{i,j}$ are coefficients in $R_q$. Thus, by
Lemma~\ref{identity},
\begin{multline*}
\sum_{x=0}^y p(x) q^x  = \sum_{x=0}^y \sum_{i=0}^{t-1}
\sum_{j=0}^d a_{ij} [x;q]_j \rho^{i x} q^x  = \sum_{i=0}^{t-1}
\sum_{j=0}^d a_{i,j} \rho^i \left. \frac{d}{d_q \sigma^j} \left(
\frac{\sigma^{j-1} ((\sigma q)^{y+1} - 1)}{(\sigma q-1)} \right)
\right|_{\sigma=\rho^i}.
\end{multline*}
The assertion follows after observing that
$$
\left. \frac{d}{d_q \sigma^j} \left( \frac{\sigma^{j-1} ((\sigma
q)^{y+1} - 1)}{(\sigma q-1)} \right) \right|_{\sigma=\rho^i}
$$
is a $q$-quasi-polynomial in $y$ of degree at most $j+1$. \qed

\medskip

Next we consider the second important property of polynomials for
our method. It suffices to derive an analog for $q$-polynomials.
Suppose $p(X)$ is a $q$-polynomial over $R$ and $a$ is an integer
zero of $p(X)$. Then there exists a $q$-polynomial $r(X)$ over $R$ with
$$
p(X)=([X;q]-[a;q]) \, r(X) = q^{a} [X-a;q] \, r(X).
$$
The proof follows from the following identity
$$
[X;q]^{n} - [a;q]^{n} = ([X;q]-[a;q]) \sum_{i=0}^{n-1} [X;q]^{i}
[a;q]^{n-1-i} = q^{a} [X-a;q] \sum_{i=0}^{n-1} [X;q]^{i}
[a;q]^{n-1-i}.$$ This property implies that for an integral domain
$R$ and distinct zeros $a_1, a_2, \ldots, a_r$ of the
$q$-polynomial $p(X)$ there exists a $q$-polynomial $r(X)$ with
$$
p(X)= \left( \prod_{i=1}^r [X - a_i;q] \right) r(X).
$$
This will be fundamental for the ``$q$-Lagrange interpolation'' we
use in Lemma~\ref{lagrange}.

\section{$F_q(n-1,n,c,p;k)$ is a $q$-quasi-polynomial in $k$}
\label{isqquasi}

In this section we show that $F_q(r,n,c,p;k_1,\ldots,k_{n-r})$ is
a $q$-quasi-polynomial in $k_1,k_2,\ldots,k_{n-r}$ with period
$2$. Moreover we show that the degree in $k_i$ is at most $2r$.

The following recursion is fundamental.
\begin{multline}
\label{rec} F_q(r,n,c,p;k_1,k_2,\ldots,k_{n-r})= \\
\sum_{l_1=0}^{k_1} \sum_{l_2=k_1}^{k_2} \sum_{l_3=k_2}^{k_3}
\ldots \sum_{l_{n-r}=k_{n-r-1}}^{k_{n-r}}
\sum_{l_{n-r+1}=k_{n-r}}^c
F_{q}(r-1,n,c,p;l_1,l_2,\ldots,l_{n-r+1}) \,
q^{l_{1}+l_{2}+\ldots+l_{n-r+1}}.
\end{multline} Moreover we have
\begin{multline*}
F_q(0,n,c,p;k_1,\ldots,k_n)=
\begin{cases}
1 & \text{if exactly $p$ of $k_1,k_2,\ldots,k_n$ are odd} \\
0 & \text{otherwise}
\end{cases}
\\
=\sum_{1 \le i_{1} < i_{2} < \ldots < i_{p} \le n} \prod_{j=1}^{p}
\frac{e_{1,2}(k_{i_{j}})}{e_{0,2}(k_{i_{j}})} \prod_{j=1}^{n}
e_{0,2}(k_{j})
 =:S(n,p)(k_1,\ldots,k_n),
\end{multline*}
where $x \to e_{i,t}(x)$ is the function defined on integers with
$$
e_{i,t}(x) =
\begin{cases}
1 & x \equiv i \mod t  \\
0 & \text{otherwise}
\end{cases}
= \prod_{0 \le j \le t-1, j \not= i} \frac{\rho^x -
\rho^j}{\rho^i-\rho^j},
$$
where $\rho \in \mathbb{C}$ is a primitive $t$-th root of unity. The identity
$$
F_q(0,n,c,p;k_1,\ldots,k_i,\ldots,k_n)=F_q(0,n,c,p;k_1,\ldots,k_{i}+2,\ldots,k_n)
$$
for all $i$, $1 \le i \le n$, implies that
$F_q(0,n,c,p;k_1,\ldots,k_n)$ is a $q$-quasi-polynomial with
period $2$. The recursion \eqref{rec} and Lemma~\ref{q-quasi-sum}
implies (inductively with respect to $r$) that $F_q(r,n,c,p;.)$ is
a $q$-quasi-polynomial in $(k_1,k_2,\ldots,k_{n-r})$ with period
at most $2$.

For our purpose it is convenient to define the following
generalization of $F_q(r,n,c,p;.)$.

\begin{defi} Let $n,r$, $r \le n$, be non-negative integers and $A(k_{1},\ldots,k_{n})$
  a function on $\mathbb{Z}^n$. We define $G_q(r,n,c,A)$ inductively
  with respect to $r$: $G_q(0,n,c,A)=A$ and
\begin{multline}
\label{rec1}
 G_q(r,n,c,A)(k_{1},\ldots,k_{n-r}) = \\
\sum_{l_{1}=0}^{k_{1}} \sum_{l_{2}=k_{1}}^{k_{2}} \ldots
\sum_{l_{n-r+1}=k_{n-r}}^{c}
G_q(r-1,n,c,A)(l_{1},l_{2},\ldots,l_{n-r+1}) \, q^{l_{1}+l_{2}+
\ldots +l_{n-r+1}}
\end{multline}
\end{defi}

With this definition we have
$$
F_q(r,n,c,p;k_1,\ldots,k_{n-r}) =
G_q(r,n,c,S(n,p))(k_1,\ldots,k_{n-r}).
$$

We define
$$
T(n,i) = \sum_{1 \le j_{1} < j_{2} < \ldots < j_{i} \le n}
(-1)^{k_{j_{1}}+k_{j_{2}}+\ldots+k_{j_{i}}}.
$$
The following lemma  shows that  $S(n,p)$ is a linear combination
of $T(n,1)$,$T(n,2),\ldots$,$T(n,n)$ and $T(n,0):=1$.

\begin{lem}
\label{linearcombination}
$$
S(n,p)=\frac{1}{2^{n}} \left( \sum_{i=0}^{n}
  \sum_{l=\max(0,i-n+p)}^{\min(p,i)} (-1)^{l} \, \binom{i}{l} \,
  \binom{n-i}{p-l} \, T(n,i) \right)
$$
\end{lem}

{\it Proof.} Set $[n]:=\{1,2,\ldots,n\}$ and fix $P \subseteq [n]$ with
$|P|=p$. Then
\begin{multline*}
\prod_{j \in P} e_{1,2}(k_{j}) \prod_{j \in [n] \setminus P}
e_{0,2}(k_{j})= \prod_{j \in P} \frac{1-(-1)^{k_{j}}}{2} \prod_{j
\in [n] \setminus P}
\frac{1+(-1)^{k_{j}}}{2} =  \\
= \frac{1}{2^{n}} \sum_{i=0}^{n}
\sum_{l=\max(0,i-n+p)}^{\min(p,i)} (-1)^{l} \sum_{1 \le j_{1} <
\ldots < j_{l} \le n, \atop  j_{x} \in P}
(-1)^{k_{j_{1}}+\ldots+k_{j_{l}}} \sum_{1 \le m_{1} < \ldots <
m_{i-l} \le n, \atop m_{x} \in [n] \setminus P}
(-1)^{k_{m_{1}}+\ldots+k_{m_{i-l}}},
\end{multline*}
where the second equation follows by expanding the product. In
the summation index $i$ counts the number of $\pm
(-1)^{k_{x}}$ we choose from the product of the $n$ factors of the form $1 \pm
(-1)^{k_{x}}$ and the index $l$ counts the number of
$-(-1)^{k_{x}}$ we choose. Observe that
\begin{multline*}
\sum_{P \subseteq [n], \atop |P|=p} \sum_{1 \le j_{1} < \ldots <
j_{l} \le n, \atop  j_{x} \in P} (-1)^{k_{j_{1}}+\ldots+k_{j_{l}}}
\sum_{1 \le m_{1} < \ldots < m_{i-l} \le n, \atop m_{x} \in [n]
\setminus P} (-1)^{k_{m_{1}}+\ldots+k_{m_{i-l}}} = \binom{i}{l}
\binom{n-i}{p-l} T(n,i),
\end{multline*}
because every $(-1)^{k_{x_{1}}+\ldots+k_{x_{i}}}$, $1 \le x_{1} <
\ldots < x_{i} \le n$, appears with multiplicity $\binom{i}{l}
\binom{n-i}{p-l}$ on the left-hand-side, since there are
$\binom{i}{l}$ ways to choose the elements from
$\{x_{1},\ldots,x_{i}\}=:I$ which lie in $P$ and
$\binom{n-i}{p-l}$ ways to choose the elements in $[n] \setminus
I$ which lie in $P$. The assertion follows. \qed

\medskip

Lemma~13 from \cite{fischer} implies that $G_q(n-1,n,c,1)(k)$ is a
$q$-polynomial  of degree $2n-2$ at most in $k$. More general we aim to show
that the degree of $G_q(n-1,n,c,T(n,p))(k)$ in $k$ is at most
$2n-2$ as well. (Thus our result reproves Lemma~13 from
\cite{fischer}.) The linearity of $A \to G_q(r,n,c,A)$ and
Lemma~\ref{linearcombination} then implies that the degree of
$G_q(n-1,n,c,S(n,p))$ is at most $2n-2$ in $k$.

In fact we show that the degree of $G_q(r,n,c,T(n,p))$ in
$k_i$ is at most $2r$. This degree estimation is rather
complicated. Assume by induction with respect to $r$ that the
degree of $G_q(r-1,n,c,T(n,p))(k_{1},\ldots,k_{n-r})$ in $k_{i}$ is at most $2r-2$ as
well as the degree in $k_{i+1}$. The degree of
$G_q(r,n,c,T(n,p))$ in $k_i$ is at most the degree of
$$
\sum_{l_{i}=k_{i-1}}^{k_{i}} \sum_{l_{i+1}=k_{i}}^{k_{i+1}}
G_q(r-1,n,c,T(n,p))(l_{1},\ldots,l_{n-r+1})
$$
in $k_{i}$ with $k_{0}=0$ and $k_{n-r+1}=c$. By Lemma~\ref{q-quasi-sum} this allows us to conclude
easily that the degree of $G_q(r,n,c,T(n,p))$ in $k_i$ is at most
$4r-2$, however, we want to establish that the degree is at most
$2r$. The following lemma is fundamental for this purpose. In
order to state it we need to define an operator $D_i$ which is
crucial for the analysis of \eqref{rec1}.

\begin{defi}
Let $G(k_{1},\ldots,k_{m})$ be a function in $m$ variables and $1
\le i \le m-1$. We set
\begin{multline*}
D_{i} G(k_{1},\ldots,k_{m})= \\
G(k_{1},\ldots,k_{i-1},k_{i},k_{i+1},k_{i+2},\ldots,k_{m}) +
G(k_{1},\ldots, k_{i-1},k_{i+1}+1,k_{i}-1,k_{i+2},\ldots,k_{m}).
\end{multline*}
\end{defi}

The following lemma shows the importance of this operator for the
degree estimation.

\begin{lem}
\label{D} Let $F(x_1,x_2)$ be a $q$-quasi-polynomial in $x_1$ and
$x_2$ which is in $x_1$ as well as in $x_2$ of degree  at most $R$.
Moreover assume that $D_1 F(x_1,x_2)$ is of degree at most $R$ as a
$q$-quasi-polynomial in $x_1$ and $x_2$, i.e. the linear
combination of ``monomials'' $[x_1;q]_{m} [x_2;q]_{n} \rho_1^{x_{1}}
\rho_2^{x_{2}}$ with $m+n \le R$ and where $\rho_1$ and $\rho_2$
are roots of unity. Then $ \sum_{x_1=a}^y \sum_{x_2=y}^b
F(x_1,x_2) q^{x_1+x_2}$ is of degree at most $R+2$ in $y$.
\end{lem}

{\it Proof.} Set $F_1(x_1,x_2)=D_1 F(x_1,x_2)/2$ and
$F_2(x_1,x_2)=(F(x_1,x_2)-F(x_2+1,x_1-1))/2$. Clearly
$F(x_1,x_2)=F_1(x_1,x_2)+F_2(x_1,x_2)$. Observe that
$F_2(x_2+1,x_1-1)=-F_2(x_1,x_2)$. Thus $F_2(x_1,x_2)$ is a linear
combination of expressions of the form
$$
[x_1;q]_m [x_2 +1;q]_n \rho_{1}^{x_{1}-1} \rho_{2}^{x_{2}}-
[x_1;q]_n [x_2+1;q]_m \rho_{1}^{x_{2}} \rho_{2}^{x_{1}-1}
$$
with $m,n \le R$ and where $\rho_1$
and $\rho_2$ are roots of unity . We set
$$
c(y,n,\rho)= \frac{d}{d_q \rho^n} \left( \frac{\rho^{n-1} ((\rho
q)^{y+1} - 1)}{(\rho q-1)} \right)
$$
Lemma~\ref{identity} implies
\begin{multline*}
\sum_{x_1=a}^y \sum_{x_2=y}^b ([x_1;q]_m [x_2+1;q]_n
\rho_1^{x_{1}-1}
\rho_2^{x_{2}} - [x_1;q]_n [x_2+1;q]_m \rho_1^{x_{2}} \rho_2^{x_{1}-1})
q^{x_1+x_2+1}  \\
= (c(y,m,\rho_1) - c(a-1,m,\rho_1))(c(b+1,n,\rho_2)-c(y,n,\rho_2)) \\
- (c(y,n,\rho_2) - c(a-1,n,\rho_2))(c(b+1,m,\rho_1)-c(y,m,\rho_1))
\\
= c(y,m,\rho_1) c(b+1,n,\rho_2) - c(a-1,m,\rho_1) c(b+1,n,\rho_2) +
c(a-1,m,\rho_1) c(y,n,\rho_2) - \\
c(y,n,\rho_2) c(b+1,m,\rho_1) + c(a-1,n,\rho_2) c(b+1,m,\rho_1) -
c(a-1,n,\rho_2) c(y,m,\rho_1).
\end{multline*}
Observe that $c(x,n,\rho_1)$ is a $q$-quasi-polynomial in $x$ of
degree at most $n+1$ and thus $$\sum_{x_1=a}^y \sum_{x_2=y}^b
F_2(x,y) q^{x+y} $$ is of degree at most $R+1$ in $y$. By the
assumption in the lemma $\sum_{x_1=a}^y \sum_{x_2=y}^b F_1(x,y)
q^{x+y} $ is of degree at most $R+2$ in $y$ and the assertion
follows. \qed

\medskip

\begin{lem}
\label{fund} Let $m$ be a positive integer, $1 \le i \le m$ and
$G({\bf l})$ be a function in ${\bf l}=(l_{1},\ldots,l_{m})$. Then
\begin{multline*}
D_i  \sum_{l_{1}=k_{1}}^{k_{2}} \sum_{l_{2}=k_{2}}^{k_{3}}
\ldots \sum_{l_{m}=k_{m}}^{k_{m+1}} G(l_{1},\ldots,l_{m}) \\
= - \frac{1}{2} \left( \sum_{l_1=k_{1}}^{k_2} \ldots
\sum_{l_{i-2}=k_{i-2}}^{k_{i-1}} \sum_{l_{i-1}=k_i+1}^{k_{i+1}+1}
\sum_{l_{i}=k_i}^{k_{i+1}} \sum_{l_{i+1}=k_{i}-1}^{k_{i+2}}
\sum_{l_{i+2}=k_{i+2}}^{k_{i+3}} \ldots
\sum_{l_{m}=k_{m}}^{k_{m+1}} D_{i-1} G ({\bf l}) \right. \\
\left. + \sum_{l_1=k_{1}}^{k_2} \ldots
\sum_{l_{i-2}=k_{i-2}}^{k_{i-1}} \sum_{l_{i-1}=k_{i-1}}^{k_{i}}
\sum_{l_{i}=k_{i}}^{k_{i+1}} \sum_{l_{i+1}=k_{i}-1}^{k_{i+1}-1}
\sum_{l_{i+2}=k_{i+2}}^{k_{i+3}} \ldots
\sum_{l_{m}=k_{m}}^{k_{m+1}} D_{i} G ({\bf l}) \right),
\end{multline*}
with $D_{0} G({\bf l})=0$ and $D_{m} G({\bf l})=0$.
\end{lem}

{\it Proof.} We set
$$
g(l_{i-1},l_{i},l_{i+1}) = \sum_{l_{1}=k_{1}}^{k_{2}} \ldots
\sum_{l_{i-2}=k_{i-2}}^{k_{i-1}} \sum_{l_{i+2}=k_{i+2}}^{k_{i+3}}
\ldots \sum_{l_m=k_m}^{k_{m+1}} G(l_{1},\ldots,l_{m}).
$$
It suffices to show the following.
\begin{multline}
\label{toshow} \sum_{l_{i-1}=k_{i-1}}^{k_{i}}
\sum_{l_{i}=k_{i}}^{k_{i+1}} \sum_{l_{i+1}=k_{i+1}}^{k_{i+2}}
g(l_{i-1},l_{i},l_{i+1}) + \sum_{l_{i-1}=k_{i-1}}^{k_{i+1}+1}
\sum_{l_{i}=k_{i+1}+1}^{k_{i}-1}
\sum_{l_{i+1}=k_{i}-1}^{k_{i+2}} g(l_{i-1},l_{i},l_{i+1}) \\
= - \frac{1}{2} \left( \sum_{l_{i-1}=k_{i}+1}^{k_{i+1}+1}
\sum_{l_{i}=k_{i}}^{k_{i+1}} \sum_{l_{i+1}=k_{i}-1}^{k_{i+2}}
g(l_{i-1},l_{i},l_{i+1}) + g(l_{i}+1,l_{i-1}-1,l_{i+1}) \right. \\
\left. + \sum_{l_{i-1}=k_{i-1}}^{k_{i}}
\sum_{l_{i}=k_{i}}^{k_{i+1}} \sum_{l_{i+1}=k_{i}-1}^{k_{i+1}-1}
g(l_{i-1},l_{i},l_{i+1}) + g(l_{i-1},l_{i+1}+1,l_{i}-1) \right)
\end{multline}
By \eqref{ext-sum} the left-hand-side of this equation is equal to
\begin{multline*}
\sum_{l_{i-1}=k_{i-1}}^{k_{i}} \sum_{l_{i}=k_{i}}^{k_{i+1}}
\sum_{l_{i+1}=k_{i+1}}^{k_{i+2}} g(l_{i-1},l_{i},l_{i+1}) -
\sum_{l_{i-1}=k_{i-1}}^{k_{i+1}+1} \sum_{l_{i}=k_{i}}^{k_{i+1}}
\sum_{l_{i+1}=k_{i}-1}^{k_{i+2}} g(l_{i-1},l_{i},l_{i+1}) \\
= \sum_{l_{i-1}=k_{i-1}}^{k_{i}} \sum_{l_{i}=k_{i}}^{k_{i+1}}
\sum_{l_{i+1}=k_{i+1}}^{k_{i+2}} g(l_{i-1},l_{i},l_{i+1}) -
\sum_{l_{i-1}=k_{i-1}}^{k_{i}} \sum_{l_{i}=k_{i}}^{k_{i+1}}
\sum_{l_{i+1}=k_{i}-1}^{k_{i+2}} g(l_{i-1},l_{i},l_{i+1}) \\
- \sum_{l_{i-1}=k_{i}+1}^{k_{i+1}+1} \sum_{l_{i}=k_{i}}^{k_{i+1}}
\sum_{l_{i+1}=k_{i}-1}^{k_{i+2}} g(l_{i-1},l_{i},l_{i+1}) \\
= - \sum_{l_{i-1}=k_{i-1}}^{k_{i}} \sum_{l_{i}=k_{i}}^{k_{i+1}}
\sum_{l_{i+1}=k_{i}-1}^{k_{i+1}-1} g(l_{i-1},l_{i},l_{i+1}) -
\sum_{l_{i-1}=k_{i}+1}^{k_{i+1}+1} \sum_{l_{i}=k_{i}}^{k_{i+1}}
\sum_{l_{i+1}=k_{i}-1}^{k_{i+2}} g(l_{i-1},l_{i},l_{i+1}).
\end{multline*}
The last expression is obviously equal to the right-hand-side of
\eqref{toshow} and the assertion of the lemma is proved. \qed

\medskip

We need another definition before we are able to prove the
key-lemma for the degree estimation.

\begin{defi} Let $r$ be a non-negative integer and $B(x,y)$ a
function in $x$ and $y$. We define $D(r,B)(x,y)$ recursively:
$D(0,B)(x,y)=B(x,y)$ and
$$
D(r,B)(x,y)= \sum_{x'=x+1}^{y+1} \sum_{y'=x}^{y} D(r-1,B)(x',y')
q^{x'+y'}.
$$
\end{defi}

The following lemma establishes a recursion which expresses $D_i
G_q(r,n,c,A)$ in terms of $G_q(r,n-2,c+2,A'_{j})$ and $D(r,B_{j})$ if $A$
fulfills a certain ``decomposition condition''.

\begin{lem}
\label{heart} Let $n,r$ be integers with $0 \le r \le n$ and
$A(k_{1},\ldots,k_{n})$ a function on $\mathbb{Z}^{n}$. Assume
that there exist two families of functions $(B_{j}(x,y))_{1 \le j
\le m}$ and $(A'_{j}(k_{1},\ldots,k_{n-2}))_{1 \le j \le m}$ with
the property that
$$
D_{i} A(k_{1},\ldots,k_{n}) =
\sum_{j=1}^{m}
B_{j}(k_{i}+i,k_{i+1}+i) A'_{j}(k_{1},\ldots,k_{i-1},k_{i+2}+2,\ldots,k_{n}+2)
$$
for all $i$, $1 \le i \le n-1$. Then
\begin{multline*}
D_{i} G(r,n,c,A)(k_{1},\ldots,k_{n-r}) =
\sum_{j=1}^{m}
\frac{(-1)^{r}}{2^{r}} q^{r(-2 n +r + 1)} \, D(r,B_{j})(k_{i}+i,k_{i+1}+i) \\
\times G(r,n-2,c+2,A'_{j})(k_{1},\ldots,k_{i-1},k_{i+2}+2,\ldots,k_{n-r}+2).
\end{multline*}
\end{lem}

{\it Proof.} We show the assertion by induction with respect to
$r$. For $r=0$ there is nothing to prove. Thus we assume $r > 0$.
By the induction hypothesis we may assume that
\begin{multline}
\label{expr1}
D_{i} G(r-1,n,c,A)(l_{1},\ldots,l_{n-r+1}) = \sum_{j=1}^{m}
\frac{(-1)^{r-1}}{2^{r-1}} q^{(r-1)(-2n+r)} D(r-1,B_{j})(l_{i}+i,l_{i+1}+i) \\
\times
G(r-1,n-2,c+2,A'_{j})(l_{1},\ldots,l_{i-1},l_{i+2}+2,\ldots,l_{n-r+1}+2)
\end{multline}
and
\begin{multline}
\label{expr2}
D_{i+1} G(r-1,n,c,A)(l_{1},\ldots,l_{n-r+1}) = \\
\sum_{j=1}^{m}
\frac{(-1)^{r-1}}{2^{r-1}} q^{(r-1)(-2n+r)} D(r-1,B_{j})(l_{i+1}+i+1,l_{i+2}+i+1) \\
\times
G(r-1,n-2,c+2,A'_{j})(l_{1},\ldots,l_{i},l_{i+3}+2,\ldots,l_{n-r+1}+2).
\end{multline}
By Lemma~\ref{fund} we have
\begin{multline*}
D_{i} G(r,n,c,A)(k_{1},\ldots,k_{n-r})  = \\
-\frac{1}{2}
 \left( \sum_{l_1=0}^{k_1} \ldots
\sum_{l_{i-1}=k_{i-2}}^{k_{i-1}} \sum_{l_{i}=k_i+1}^{k_{i+1}+1}
\sum_{l_{i+1}=k_i}^{k_{i+1}} \sum_{l_{i+2}=k_{i}-1}^{k_{i+2}}
\ldots \sum_{l_{n-r+1}=k_{n-r}}^{c} D_{i} G(r-1,n,c,A)({\bf l}) q^{l_{1}+\ldots+l_{n-r+1}} \right. \\
\left. + \sum_{l_1=0}^{k_1} \ldots
\sum_{l_{i}=k_{i-1}}^{k_{i}} \sum_{l_{i+1}=k_{i}}^{k_{i+1}}
\sum_{l_{i+2}=k_{i}-1}^{k_{i+1}-1}
\sum_{l_{i+3}=k_{i+2}}^{k_{i+3}} \ldots
\sum_{l_{n-r+1}=k_{n-r}}^{c} D_{i+1} G(r-1,n,c,A) ({\bf l}) q^{l_{1}+\ldots+l_{n-r+1}}\right)
\end{multline*}
In this expression we replace  $D_{i} G(r-1,n,c,A)({\bf l})$ by \eqref{expr1}
and  $D_{i+1} G(r-1,n,c,A) ({\bf l})$ by \eqref{expr2}, and obtain
\begin{multline*}
\frac{(-1)^{r}}{2^{r}} q^{(r-1)(-2n+r)} \sum_{j=1}^{m}
\left(\sum_{l_{i}=k_{i}+1}^{k_{i+1}+1} \sum_{l_{i+1}=k_{i}}^{k_{i+1}}
    D(r-1,B_{j})(l_{i}+i,l_{i+1}+i) q^{l_{i}+l_{i+1}} \right. \\
\times \sum_{l_{1}=0}^{k_{1}} \ldots \sum_{l_{i-1}=k_{i-2}}^{k_{i-1}}
\sum_{l_{i+2}=k_{i}-1}^{k_{i+2}} \ldots \sum_{l_{n-r+1}=k_{n-r}}^{c} \\
G(r-1,n-2,c+2,A'_{j})(l_{1},\ldots,l_{i-1},l_{i+2}+2,\ldots,l_{n-r+1}+2)
q^{l_{1}+\ldots+l_{i-1}+l_{i+2}+\ldots+l_{n-r+1}} \\
+
\sum_{l_{i+1}=k_{i}}^{k_{i+1}} \sum_{l_{i+2}=k_{i}-1}^{k_{i+1}-1}
D(r-1,B_{j})(l_{i+1}+i+1,l_{i+2}+i+1) q^{l_{i+1}+l_{i+2}} \\
\times
\sum_{l_{1}=0}^{k_{1}} \ldots \sum_{l_{i}=k_{i-1}}^{k_{i}}
\sum_{l_{i+3}=k_{i+2}}^{k_{i+3}} \ldots \sum_{l_{n-r+1}=k_{n-r}}^{c} \\
G(r-1,n-2,c+2,A'_{j})(l_{1},\ldots,l_{i},l_{i+3}+2,\ldots,l_{n-r+1}+2)
q^{l_{1}+\ldots+l_{i}+l_{i+3}+\ldots+l_{n-r+1}} \Bigg) \\
=
\frac{(-1)^{r}}{2^{r}} q^{r(-2n+r+1)} \sum_{j=1}^{m}
D(r,B_{j})(k_{i}+i,k_{i+1}+i) \\
\times \Bigg(
\sum_{l_{1}=0}^{k_{1}} \ldots \sum_{l_{i-1}=k_{i-2}}^{k_{i-1}}
\sum_{l_{i+2}=k_{i}+1}^{k_{i+2}+2} \sum_{l_{i+3}=k_{i+2}+2}^{k_{i+3}+2}\ldots \sum_{l_{n-r+1}=k_{n-r}+2}^{c+2} \\
G(r-1,n-2,c+2,A'_{j})(l_{1},\ldots,l_{i-1},l_{i+2},\ldots,l_{n-r+1})
q^{l_{1}+\ldots+l_{i-1}+l_{i+2}+\ldots+l_{n-r+1}}  \\
+
\sum_{l_{1}=0}^{k_{1}} \ldots \sum_{l_{i-1}=k_{i-2}}^{k_{i-1}} \sum_{l_{i}=k_{i-1}}^{k_{i}}
\sum_{l_{i+3}=k_{i+2}+2}^{k_{i+3}+2} \ldots \sum_{l_{n-r+1}=k_{n-r}+2}^{c+2}
\\
G(r-1,n-2,c+2,A'_{j})(l_{1},\ldots,l_{i},l_{i+3},\ldots,l_{n-r+1})
q^{l_{1}+\ldots+l_{i}+l_{i+3}+\ldots+l_{n-r+1}} \Bigg)
\end{multline*}
The last expression is obviously equal to the right-hand-side of
the equation in the statement of the lemma. \qed

\medskip

In the next lemma we give a bound for the degree of
$D(r,B)(x,y)$.

\begin{lem}
\label{degree} Suppose $B(x,y)$ is a $q$-quasi-polynomial in $x$
and $y$ of degree $d$, i.e. the linear combination of terms of the form $[x;q]_m
[y;q]_n \rho_1^x \rho_2^y$ with $m+n \le d$ and $\rho_1, \rho_2$
are roots of unity. Then $D(r,B)(x,y)$ is of degree at most $2r+d$
in $x$ and $y$.
\end{lem}

{\it Proof.} By induction with respect to $r$ it suffices to show
that
$$
\sum_{x'=x+1}^{y+1} \sum_{y'=x}^{y} [x';q]_m [y';q]_n
\rho_1^{x'-1} \rho_2^{y'-1} q^{x'+y'}
$$
is of degree at most $m+n+2$ in $x$ and $y$. Using the notation from
Lemma~\ref{q-quasi-sum} we see that this double sum is equal to
\begin{eqnarray*}
\left( \sum_{x'=x+1}^{y+1} [x';q]_{m} \rho_{1}^{x'-1} q^{x'}
\right) \left( \sum_{y'=x}^{y} [y';q]_{n} \rho_{2}^{y'-1} q^{y'}
\right) \\
= (c(y+1,m,\rho_{1}) - c(x,m,\rho_{1}))(c(y,n,\rho_{2}) -
c(x-1,n,\rho_{2}))
\end{eqnarray*}
and the assertion follows. \qed

\medskip

In order to apply Lemma~\ref{heart} to our situation we show that
$T(n,p)(k_{1},\ldots,k_{n})$ has the ``decomposition property'' from Lemma~\ref{heart}. Observe that
\begin{multline*}
D_{i} T(n,p) = 2 \, T(n-2,p)(k_{1},\ldots,k_{i-1},k_{i+2}+2,\ldots,k_{n}+2) \\ +
(-1)^{k_{i}+k_{i+1}} 2 \,
T(n-2,p-2)(k_{1},\ldots,k_{i-1},k_{i+2}+2,\ldots,k_{n}+2),
\end{multline*}
where $T(n,p)=0$ if $p < 0$ or $p > n$.
Thus, by Lemma~\ref{heart},
\begin{multline*}
D_{i} G_q(r,n,c,T(n,p)) = \frac{(-1)^r}{2^r} q^{r(-2n+r+1)} \big(
D(r,1)(k_i+i,k_{i+1}+i) \\ \times G_q(r,n-2,c+2,2 \,
T(n-2,p))(k_1,\ldots,k_{i-1},k_{i+2}+2,\ldots,k_{n-r}+2)
\\  + D(r,(-1)^{k_i+k_{i+1}})(k_i+i,k_{i+1}+i) \\ \times G_q(r,n-2,c+2,2 \,
T(n-2,p-2))(k_1,\ldots,k_{i-1},k_{i+2}+2,\ldots,k_{n-r}+2) \big).
\end{multline*}
By Lemma~\ref{degree} $D(r,1)(k_{i},k_{i+1})$ as well as
$D(r,(-1)^{k_i+k_{i+1}})(k_{i},k_{i+1})$ are $q$-quasi-polynomials
in $k_i$ and $k_{i+1}$ of degree at most $2 r$ and thus the same
is true for $D_{i} G_q(r,n,c,T(n,p))$. Finally we show that this
implies that $G_q(r,n,c,T(n,p))$ is a  $q$-quasi-polynomial in
$k_i$ of degree at most $2 r$ for all $i$.

\begin{lem}
\label{degree1} Let $n,r$ be positive integers, $r \le n$ and $0 \le p \le
n$. Then $G_q(r,n,c,T(n,p))$ is a $q$-quasi-polynomial in $k_i$ of
degree at most $2r$ for $i=1,2,\ldots,n-r$.
\end{lem}

{\it Proof.}  We show the assertion by induction with respect to
$r$. For $r=0$ there is nothing to prove. We assume that $r > 0$
and that the assertion is true for $G_q(r-1,n,c,T(n,p))$. The degree of
$D_i G_q(r-1,n,c,T(n,p))(l_1,\ldots,l_{n-r+1})$ as a $q$-quasi-polynomial
in $l_i$ and $l_{i+1}$ is at most $2r-2$. Therefore, by
Lemma~\ref{D}, the degree of
$$
\sum_{l_i=k_{i-1}}^{k_i} \sum_{l_{i+1}=k_i}^{k_{i+1}}
G_q(r-1,n,c,T(n,p))(l_1,\ldots,l_{n-r+1}) q^{l_i+l_{i+1}}
$$
in $k_i$ is at most $2r$. By \eqref{rec1} the same is true for the
degree of $G_q(r,n,c,T(n,p))$ in $k_i$. \qed

\begin{cor}
\label{degree2}
Let $n$ be a positive integer and $0 \le p \le n$.
Then $F_q(n-1,n,c,p;k)$ is a $q$-quasi-polynomial over $\mathbb{C}$
of degree at most $2n -2$ in $k$.
\end{cor}

\section{Exploring properties of the $q$-quasi-polynomial $F_q(n-1,n,c,p;k)$}
\label{properties}

First we observe that $F_q(n-1,n,c,p;k)$ is zero for
$k=-1,-2,\ldots,-n+1$ and $k=c+1,c+2,\ldots,c+n-1$.

\begin{lem}
\label{zeros} Let $r,n,p$ be integers, $0 \le r \le n$ and $0 \le
p \le n$. Then $F_q(r,n,c,p;.)$ is zero for $k_1=-1,-2,\ldots,-r$
and $k_{n-r}=c+1,c+2,\ldots,c+r$.
\end{lem}

{\it Proof.} It suffices to show that there exists no
$(r,n,c)$-pattern with first row
$$(0,k_1,\ldots,k_{n-r},c),$$
if $k_1=-1,-2,\ldots,-r$ or $k_{n-r}=c+1,c+2,\ldots,c+r$. Indeed,
suppose $(a_{i,j})$ is an $(r,n,c)$-pattern with $a_{r+1,r+1} \in
\{-1,-2,\ldots,-r\}$. In particular we have $0 > a_{r+1,r+1}$ and
thus the definition of $(r,n,c)$-patterns implies that $0 >
a_{r,r} > a_{r+1,r+1}$. In a similar way we obtain $0
> a_{1,1} > a_{2,2} > \ldots > a_{r,r}
> a_{r+1,r+1}$. This is, however, a contradiction, since there
exist no $r$ distinct integers between $0$ and $a_{r+1,r+1}$. The
case that $a_{r+1,n} \in \{c+1,c+2,\ldots,c+r\}$ is similar. \qed

\medskip

The zeros in Lemma~\ref{zeros} do not determine the
$q$-quasi-polynomial $F_q(n-1,n,c,p;k)$ uniquely and thus we need
additional properties. To this end we have the following two
lemmas.

\begin{lem}
\label{DT} Let $r$ be a non-negative integer. Then we have
$$
D(r,(-1)^{x+y}) = (-1)^{x+y} q^{r(x+y)} C + T(x,y),
$$
where $C \in \mathbb{Q}(q)$ and $T(x,y)$ is a $q$-polynomial in
$x$ and $y$ over $\mathbb{Q}$.
\end{lem}

{\it Proof.} The assertion follows from the following identity by
induction with respect to $r$.
$$
\sum_{x'=x+1}^{y+1} \sum_{y'=x}^{y} (-1)^{x'+y'} Q^{x'+y'} =
\frac{-Q^{2x+1} - Q^{2y+3} - 2 (-1)^{x+y} Q^{x+y+2}}{(1+Q)^{2}} \qed
$$

\medskip

Suppose $p(X)$ is a $q$-quasi-polynomial in $X$ with period $1$ or $2$. Then there
exist unique $q$-polynomials $p_{1}(X)$ and $p_{2}(X)$ with the property
that
$$
p(X) = (-1)^{X} p_{1}(X) + p_{2}(X).
$$
We say that $p_{1}(X)$ is the signed part of $p(X)$, in symbols
$\SP_{X} (p(X))=p_{1}(X)$. The following lemma shows that the
signed part of the $q$-quasi-polynomials $F_{q}(r,n,c,p;.)$ have a
quite simple structure.

\begin{lem} Let $r,n,i,p$ be integers, $r$ non-negative, $n$ positive, $1 \le
  i < n-r$ and $0 \le p \le n$. Then
$$ \SP_{k_{i}} F_{q}(r,n,c,p;k_{1},\ldots,k_{n-r}) / q^{r k_{i}}$$
is independent of $k_{i}$.
\end{lem}

{\it Proof.} We show the assertion by induction with respect to $r$. For
$r=0$ there is nothing to prove. Let $r > 0$. It suffices to prove that
\begin{equation}
\label{SP1} \SP_{k_{i}} \left( \sum_{l_{i}=k_{i-1}}^{k_{i}}
\sum_{l_{i+1}=k_{i}}^{k_{i+1}} D_{i}
G(r-1,n,c,T(n,p))(l_{1},\ldots,l_{n-r+1}) q^{l_{i}+l_{i+1}}
\right) /q^{r k_{i}}
\end{equation}
and
\begin{multline}
\label{SP2}
\SP_{k_{i}} \Big( \sum_{l_{i}=k_{i-1}}^{k_{i}} \sum_{l_{i+1}=k_{i}}^{k_{i+1}} \Big(
G(r-1,n,c,T(n,p))(l_{1},\ldots,l_{i},l_{i+1},\ldots,l_{n-r+1})  \\
\ -
G(r-1,n,c,T(n,p))(l_{1},\ldots,l_{i+1}+1,l_{i}-1,\ldots,l_{n-r+1})
\Big) q^{l_{i}+l_{i+1}}  \Big) /q^{r k_{i}}
\end{multline}
are independent of $k_{i}$, where $k_{0}=0$ and $k_{n-r+1}=c$.
By Lemma~\ref{heart} and Lemma~\ref{DT} $D_{i}
G(r-1,n,c,T(n,p))(l_{1},\ldots,l_{n-r+1})$ is of the
form
$$(-1)^{l_{i}+l_{i+1}} q^{(r-1)(l_{i}+l_{i+1})}
H(l_{1},\ldots,l_{i-1},l_{i+2},\ldots,l_{n-r+1})  +
T(l_{1},\ldots,l_{n-r+1}),$$ where
$H(l_{1},\ldots,l_{i-1},l_{i+2},\ldots,l_{n-r+1})$ is a
$q$-quasi-polynomial, $T(l_{1},\ldots,l_{n-r+1})$ is a
$q$-quasi-polynomial in
$(l_{1},\ldots,l_{i-1},l_{i+2},\ldots,l_{n-r+1})$ and a
$q$-polynomial in $l_{i}$ and $l_{i+1}$. Therefore $T$ does not
contribute to the signed part of \eqref{SP1}. Moreover we have
\begin{multline*}
\sum_{l_{i}=k_{i-1}}^{k_{i}} \sum_{l_{i+1}=k_{i}}^{k_{i+1}}
(-1)^{l_{i}+l_{i+1}} q^{(r-1)(l_{i}+l_{i+1})}
H(l_{1},\ldots,l_{i-1},l_{i+2},\ldots,l_{n-r+1}) \, q^{l_{i}+l_{i+1}} \\
=\frac{1}{(1+q^r)^2} ((-1)^{k_{i}} q^{r \, k_{i}} (-1)^{k_{i-1}}
q^{r \, k_{i-1}} +
 (-1)^{k_{i}} q^{r k_{i}} (-1)^{k_{i+1}} q^{r (k_{i+1}+2)} \\
+ q^{r (2 k_{i} + 1)} + (-1)^{k_{i-1}+k_{i+1}} q^{r \, (k_{i-1} + k_{i+1}
  +1)}) \,
H(l_{1},\ldots,l_{i-1},l_{i+2},\ldots,l_{n-r+1})
\end{multline*}
and the first assertion follows. For the second assertion observe that by the
induction hypothesis $G(r-1,n,c,T(n,p))(l_{1},\ldots,l_{n-r+1})$ is a linear
combination of expressions of the form
$$
[l_{i};q]_{m} [l_{i+1}+1;q]_{n},
$$
$$
q^{(r-1) l_{i}} (-1)^{l_{i}-1} [l_{i+1}+1;q]_{n}
$$
and
$$
q^{(r-1) l_{i}} (-1)^{l_{i}-1} q^{(r-1)(l_{i+1}+1)} (-1)^{l_{i+1}}
$$
over $R_{qq}[l_{1},\ldots,l_{i-1},l_{i+2},\ldots,l_{n-r}]$
Therefore
\begin{multline}
G(r-1,n,c,T(n,p))(l_{1},\ldots,l_{i},l_{i+1},\ldots, l_{n-r+1}) \\ -
G(r-1,n,c,T(n,p))(l_{1},\ldots,l_{i+1}+1,l_{i}-1,\ldots,
l_{n-r+1})
\end{multline}
is a linear combination of expressions of the form
\begin{equation}
\label{second}
[l_{i};q]_{m} [l_{i+1}+1;q]_{n}
-
[l_{i};q]_{n} [l_{i+1}+1;q]_{m},
\end{equation}
\begin{equation}
\label{first}
q^{(r-1) l_{i}} (-1)^{l_{i}-1} [l_{i+1}+1;q]_{n}
- [l_{i};q]_{n} q^{(r-1)(l_{i+1}+1)} (-1)^{l_{i+1}}
\end{equation}
and
\begin{equation}
\label{third}
q^{(r-1) l_{i}} (-1)^{l_{i}-1} q^{(r-1)(l_{i+1}+1)} (-1)^{l_{i+1}} -
q^{(r-1) l_{i}} (-1)^{l_{i}-1} q^{(r-1)(l_{i+1}+1)} (-1)^{l_{i+1}}=0
\end{equation}
over $R_{qq}[l_{1},\ldots,l_{i-1},l_{i+2},\ldots,l_{n-r}]$.
Expressions of the form \eqref{second} do not contribute to the
signed part of \eqref{SP2}. For expressions of the
form~\ref{first} observe that
\begin{multline*}
\sum_{l_{i}=k_{i-1}}^{k_{i}} \sum_{l_{i+1}=k_{i}}^{k_{i+1}}
\left( q^{(r-1) l_{i}} (-1)^{l_{i}-1} [l_{i+1}+1;q]_{n}
- [l_{i};q]_{n}  q^{(r-1)(l_{i+1}+1)} (-1)^{l_{i+1}} \right)
q^{l_{i}+l_{i+1}+1} \\
= - \frac{(-1)^{k_{i-1}} q^{r \, k_{i-1}} + q^{r} (-1)^{k_{i}} q^{r k_{i}}}{1
  + q^{r}} \frac{q}{[n+1;q]} ([k_{i+1}+1;q]_{n+1} - [k_{i};q]_{n+1}) \\
- \frac{q}{[n+1;q]} ([k_{i};q]_{n+1} - [k_{i-1}-1;q]_{n+1})
\frac{q^{r} (-1)^{k_{i}} q^{r k_{i}} + q^{2 r} (-1)^{k_{i+1}} q^{r
  k_{i+1}}}{1 + q^{r}} \\
= \frac{q}{(1+q^{r}) [n+1;q]}
\big( - (-1)^{k_{i-1}} q^{r k_{i-1}} [k_{i+1}+1;q]_{n+1} -
q^{r} (-1)^{k_{i}} q^{r k_{i}} [k_{i+1}+1;q]_{n+1} \\ +
(-1)^{k_{i-1}} q^{r k_{i-1}} [k_{i};q]_{n+1}  +
[k_{i-1}-1;q]_{n+1} q^{r} (-1)^{k_{i}} q^{r k_{i}} \\ -
[k_{i};q]_{n+1} q^{2 r} (-1)^{k_{i+1}} q^{r k_{i+1}} +
[k_{i-1}-1;q]_{n+1} q^{2r} (-1)^{k_{i+1}} q^{r k_{i+1}} \big).
\end{multline*}
The assertion follows. \qed

\medskip

\begin{cor}
\label{quasi-form}
Let $n$ be a positive integer. Then
$$
F(n-1,n,c,p;k) = P_{n,c,p}(k) + (-1)^{k} q^{(n-1) k} L_{n,c,p},
$$
where $P_{n,c,p}(k)$ is a $q$-polynomial in $k$ and $L_{n,c,p}$ is
independent of $k$.
\end{cor}

\medskip

We define
$$
G_{n,c,p}=\sum_{k=0}^{c} F_{q}(n-1,n,c,p;k) q^{k}.
$$
This is the generating function of strict plane partitions with
parts in $\{1,2,\ldots,n\}$, at most $c$ columns and $p$ odd rows.
In the final lemma of this section we prove that some special
evaluations of $F_{q}(n-1,n,c,p;k)$ in $k$ can be expressed in
terms of the generating function $G_{n-1,c,p}$. This lemma,
together with Lemma~\ref{zeros} and Corollary~\ref{quasi-form},
provides enough properties in order to compute
$F_{q}(n-1,n,c,p;k)$ in the following section.

\begin{lem}
\label{initial-conditions} Let $n$ be a positive integer. If $p
\not=n$ then
$$
F_{q}(n-1,n,c,p;0) = G_{n-1,c,p}
$$
and if $p \not= 0$
$$
F_{q}(n-1,n,c,p;-n) = (-1)^{n-1} q^{-3 n (n-1)/2} G_{n-1,c+2,p-1}.
$$
Moreover we have
$$
F_{q}(n-1,n,c,n;1)=q^{(n+2)(n-1)/2}  G_{n-1,c-1,0}
$$
and
$$
F_{q}(n-1,n,c,0;-n-1)=(-1)^{n-1} q^{-(n-1)(2n+1)} G_{n-1,c+3,n-1}.
$$
\end{lem}

{\it Proof.} First let $(a_{i,j})$ be an $(n-1,n,c)$-pattern with
$a_{n,n}=0$ and exactly $p$ numbers of $a_{1,1}, a_{1,2},\ldots,a_{1,n}$
are odd. This implies that $a_{i,i}=0$ for all $i$ and thus
$(a_{i,j})_{1 \le i \le n-1, i \le j \le n+1}$ is an
$(n-2,n-1,c)$-pattern with $p$ of $a_{1,2},
a_{1,3},\ldots,a_{1,n}$ are odd. In fact this induces a norm-preserving and
sign-preserving bijection between
these $(n-1,n,c)$-patterns and these $(n-2,n-1,c)$-patterns.
The first identity is proved.

Next let $(a_{i,j})$ be an $(n-1,n,c)$-pattern with $a_{n,n}=-n$
and exactly $p$ of $a_{1,1}, a_{1,2},\ldots,a_{1,n}$ are odd. This
implies that $a_{i,i}=-i$. Therefore $a_{i,i+1} \notin
\{-3,-4,\ldots,-n\}$ for $i=1,\ldots,n-1$. If we set
$b_{i,j}:=a_{i,j}+2$ for $i < j$ and $b_{i,i}=0$ then
$(b_{i,j})_{1 \le i \le n-1, i \le j \le n+1}$ is an
$(n-2,n-1,c+2)$-pattern with $p-1$ of $b_{1,2}, b_{1,3}, \ldots,
b_{1,n}$ are odd. Again this induces a bijection. However, the
bijection is neither norm-preserving nor sign-preserving and but
the factor $(-1)^{n-1} q^{-3n(n-1)/2}$ takes into account  the
changes of norm and sign.

For the third identitiy let $(a_{i,j})$ be an $(n-1,n,c)$-pattern with
$a_{n,n}=1$ and all $a_{1,1},a_{1,2},\ldots,a_{1,n}$ are
odd. The first assumption implies that $a_{i,i} \in \{0,1\}$, the second
that $a_{1,1}=1$ and therefore $a_{i,i}=1$ for all $i$. If
we set $b_{i,j}=a_{i,j}-1$ then  $(b_{i,j})_{1 \le i \le n-1, i \le j
\le n+1}$ is an $(n-2,n-1,c-1)$-pattern, where all
$b_{1,2}, b_{1,3}, \ldots, b_{1,n}$ are even and the identity follows.

The proof of the fourth identity is similar. \qed

\section{Computation of $F_{q}(n-1,n,c,p;k)$}
\label{compute}

In this section we compute $F_{q}(n-1,n,c,p;k)$ using the
properties we have established in the previous section. For these
computations we need some $q$-summation formulas which we derive
in Section~\ref{q-sums}. First we see that
Corollary~\ref{degree2}, Lemma~\ref{zeros} and
Corollary~\ref{quasi-form} imply a first strong assertion on the
form of $F_{q}(n-1,n,c,p;k)$.

\begin{lem}
\label{lagrange}
Let $n$ be a positive integer and $0 \le p \le
n$. Then
  $F_{q}(n-1,n,c,p;k)$ is of the form
\begin{multline*}
 M_{n,c,p} [k+1;q]_{n-1} [k-c-n+1;q]_{n-1}  +
L_{n,c,p} \Bigg ((-1)^{k} q^{(n-1) k} + (-1)^n q^{(n-1)(2c+n)/2} \\
\times \sum_{i=1}^{n-1}  \Big( (-1)^{c} q^{\binom{i}{2}}
\frac{[k+1;q]_{n-1}
  [k-c-i+1;q]_{i-1} [k-c-n+1;q]_{n-i-1}}{[1;q]_{i-1} [1;q]_{n-1-i}
  [c+i+1;q]_{n-1}} \\
- q^{\binom{i}{2}} \frac{[k+1;q]_{i-1}
  [k+i+1;q]_{n-i-1} [k-c-n+1;q]_{n-1}}{ [1;q]_{i-1}
  [1;q]_{n-1-i} [c+i+1;q]_{n-1}} \Big)   \Bigg)
\end{multline*}
\end{lem}

{\it Proof.} By Lemma~\ref{zeros} and
Corollary~\ref{quasi-form} we know that for
$k \in \{-1,-2,\ldots,-n+1\}$ and $k \in \{c+1,c+2,\ldots,c+n-1\}$ we have
$$
P_{n,c,p}(k) = (-1)^{k+1} q^{(n-1) k } L_{n,c,p}.
$$
By Corollary~\ref{degree2} $P_{n,c,p}(k)$ is a $q$-polynomial in
$k$ of degree at most $2n-2$. By $q$-Lagrange interpolation the
following polynomial is the unique $q$-polynomial of degree at
most $2n-3$ with the same evaluations as $P_{n,c,p}(k)$ at $k \in
\{-1,-2,\ldots,-n+1\}$ and $k \in \{c+1,c+2,\ldots,c+n-1\}$.
\begin{multline*}
\sum_{i=-n+1}^{-1} (-1)^{i+1}  \, q^{(n-1) \, i} L_{n,c,p} \prod_{-n+1 \le j \le
  -1, j \not= i} \frac{[k-j;q]}{[i-j;q]}  \prod_{j=1}^{n-1}
\frac{[k-c-j;q]}{[i-c-j;q]} \\
+ \sum_{i=1}^{n-1} (-1)^{c+i+1}  \, q^{(n-1)(c+i)} L_{n,c,p}
\prod_{j=-n+1}^{-1} \frac{[k-j;q]}{[c+i-j;q]} \prod_{1 \le j \le
n-1, j \not= i}
\frac{[k-c-j;q]}{[i-j;q]}.
\end{multline*}
This is equal to
\begin{multline*}
\sum_{i=1}^{n-1} \Big( (-1)^{i+1}  \, q^{-(n-1) i} L_{n,c,p}
\prod_{j=1}^{i-1} \frac{[k+j;q]}{[j-i;q]} \prod_{j=i+1}^{n-1}
\frac{[k+j;q]}{[j-i;q]} \prod_{j=1}^{n-1}
\frac{[k-c-j;q]}{[-i-c-j;q]} \\
+ (-1)^{c+i+1} \, q^{(n-1)(c+i)} L_{n,c,p} \prod_{j=1}^{n-1}
\frac{[k+j;q]}{[c+i+j;q]} \prod_{j=1}^{i-1}
\frac{[k-c-j;q]}{[i-j;q]} \prod_{j=i+1}^{n-1} \frac{[k-c-j;q]}{[i-j;q]} \Big) \\
= L_{n,c,p} \sum_{i=1}^{n-1} \Big( (-1)^{i+1} q^{-(n-1) i}
\frac{[k+1;q]_{i-1} [k+i+1;q]_{n-1-i}
[k-c-n+1;q]_{n-1}}{[1-i;q]_{i-1} [1;q]_{n-1-i} [-n+1-c-i;q]_{n-1}}
\\ + (-1)^{c+i+1} q^{(n-1)(c+i)}
\frac{[k+1;q]_{n-1} [k-c-i+1;q]_{i-1}
[k-c-n+1;q]_{n-1-i}}{[c+i+1;q]_{n-1} [1;q]_{i-1}
[i-n+1;q]_{n-1-i}} \Big).
\end{multline*}
The difference of $P_{n,c,p}(k)$ and the $q$-polynomial displayed above
is a $q$-polynomial of degree $2n-2$ at most which vanishes for $k \in \{-1,\ldots,-n+1\}$ and
$k \in \{c+1,\ldots,c+n-1\}$. Thus this difference is equal to
$M_{n,c,p} [k+1;q]_{n-1} [k-c-n+1;q]_{n-1}$, where $M_{n,c,p}$ is
a factor independent of $k$, which still has to be determined.
We use the identity
$$
[z;q]_n = [-z-n+1;q]_n (-1)^n q^{n(z+(n-1)/2)}
$$
in order to obtain the expression for $P_{n,c,p}(k)$ in
the statement of the lemma. \qed

\medskip

We set
\begin{multline*}
U_{n,c}(k)= \\
(-1)^n q^{(n-1)(2c+n)/2} \sum_{i=1}^{n-1}  \Big( (-1)^{c} q^{\binom{i}{2}}
\frac{[k+1;q]_{n-1}
  [k-c-i+1;q]_{i-1} [k-c-n+1;q]_{n-i-1}}{[1;q]_{i-1} [1;q]_{n-1-i}
  [c+i+1;q]_{n-1}} \\
- q^{\binom{i}{2}} \frac{[k+1;q]_{i-1}
  [k+i+1;q]_{n-i-1} [k-c-n+1;q]_{n-1}}{ [1;q]_{i-1}
  [1;q]_{n-1-i} [c+i+1;q]_{n-1}}  \Big) + (-1)^k q^{(n-1) k}
\end{multline*}
and $W_{n,c}(k)=[k+1;q]_{n-1} [k-c-n+1;q]_{n-1}$. With these definitions
Lemma~\ref{lagrange} states that
\begin{equation}
\label{UW}
F_q(n-1,n,c,p;k)=L_{n,c,p} U_{n,c}(k) + M_{n,c,p} W_{n,c}(k).
\end{equation}
It remains to compute $L_{n,c,p}$ and $M_{n,c,p}$. In the
following lemma we give recursive formulas for $L_{n,c,p}$ and
$M_{n,c,p}$ with respect to $n$. It is an immediate consequence of
Lemma~\ref{initial-conditions}.

\begin{lem}
\label{recursion}
The initial conditions are
$L_{1,c,p}=\frac{(-1)^p}{2}$ and $M_{1,c,p}=\frac{1}{2}$. If $p
\not=0, n$ we have
$$
L_{n,c,p} = \frac{G_{n-1,c,p} W_{n,c}(-n)  + (-1)^{n} q^{-3 n (n-1)/2}
  G_{n-1,c+2,p-1} W_{n,c}(0)}{U_{n,c}(0)
W_{n,c}(-n) - U_{n,c} (-n) W_{n,c}(0)}.
$$
In case that $p=0$ we have the following recursion
$$
L_{n,c,0}= \frac{(-1)^{n-1} q^{-(n-1)(2n+1)} G_{n-1,c+3,n-1}
W_{n,c}(0) - G_{n-1,c,0} W_{n,c} (-n-1)}{U_{n,c}(-n-1) W_{n,c}(0)
- U_{n,c}(0) W_{n,c}(-n-1)},
$$
and for $p=n$ we have
$$
L_{n,c,n} = \frac{q^{(n+2)(n-1)/2} G_{n-1,c-1,0} W_{n,c}(-n) +
(-1)^{n} q^{-3 n (n-1)/2} G_{n-1,c+2,n-1} W_{n,c}(1)}
{U_{n,c}(1) W_{n,c}(-n) - U_{n,c}(-n) W_{n,c}(1)}.
$$
Concerning $M_{n,c,p}$ we have
$$
M_{n,c,p} = \frac{G_{n-1,c,p} - U_{n,c}(0) L_{n,c,p}}{W_{n,c}(0)},
$$
for $p \not= n$ and
$$
M_{n,c,p} = \frac{(-1)^{n-1} q^{- 3 n (n-1)/2} G_{n-1,c+2,p-1} -
U_{n,c}(-n) L_{n,c,p}}{W_{n,c}(-n)}
$$
for $p \not= 0$.
\end{lem}

{\it Proof.} Lemma~\ref{initial-conditions} and \eqref{UW} implies the following
equations. If $p \not= 0$ then
$$
L_{n,c,p} U_{n,c}(-n) + M_{n,c,p} W_{n,c}(-n) = (-1)^{n-1} q^{-3n
(n-1)/2} G_{n-1,c+2,p-1},
$$
and if $p \not= n$ then
$$
L_{n,c,p} U_{n,c}(0) + M_{n,c,p} W_{n,c}(0) = G_{n-1,c,p}.
$$
For $p=0$ we have
$$
L_{n,c,0} U_{n,c}(-n-1) + M_{n,c,0} W_{n,c}(-n-1) = (-1)^{n-1}
q^{-(n-1)(2n+1)} G_{n-1,c+3,n-1}
$$
and for $p = n$ we have
$$
L_{n,c,n} U_{n,c}(1) + M_{n,c,n} W_{n,c}(1) = q^{(n+2)(n-1)/2}G_{n-1,c-1,0}.
$$
For every $p \in \{0,1,2,\ldots,n\}$ this gives a system of two
linearly independent equations for $L_{n,c,p}$ and $M_{n,c,p}$. With
the help of Cramer's rule we obtain the recursions of $L_{n,c,p}$.
The recursions for $M_{n,c,p}$ are immediate consequences of the
equations. \qed

\medskip

In the following lemma we see that the denominators in the
recursive formulas for $L_{n,c,p}$ in Lemma~\ref{recursion} are
products.

\begin{lem}
\label{denominator} We have
\begin{equation*}
U_{n,c}(0) W_{n,c} (-n) - U_{n,c}(-n) W_{n,c}(0) = \frac{2
[1;q^2]_{n-1} (1+q)^{2n-1}}{q^{c(n-1)+2n(n-1)}}
\begin{cases}
\frac{[(c+2)/2;q^2]_{n-1}}{1+q} & \text{if $c$ is even} \\
\frac{[(c+1)/2;q^2]_n}{[c+n;q]} & \text{if $c$ is odd}
\end{cases},
\end{equation*}
\begin{multline*}
U_{n,c}(-n-1) W_{n,c}(0) - U_{n,c}(0) W_{n,c}(-n-1) = - \frac{ 2
[1;q^2]_{n-1} [n-1;q] (1+q)^{2n-1}}{q^{(n-1)c + 2(n-1)(n+1)}} \\
\times
\begin{cases}
[(c+2)/2;q^2]_n & \text{if $c$ is even} \\
\frac{[(c+1)/2;q^2]_{n+1} (1+q)}{[c+n+1;q]} & \text{if $c$ is odd}
\end{cases}
\end{multline*}
and
\begin{multline*}
U_{n,c}(1) W_{n,c}(-n) - U_{n,c}(-n) W_{n,c}(1)= \frac{2
[1;q^2]_{n-1} [n-1;q] (1+q)^{2n-1}}{q^{(n-1) c + n (2n -3)}
[c+n;q]} \\ \times
\begin{cases}
[c/2;q^2]_n & \text{if $c$ is even} \\
\frac{[(c-1)/2;q^2]_{n+1} (1+q)}{[c+n-1;q]} & \text{if $c$ is odd}
\end{cases}.
\end{multline*}
\end{lem}

{\it Proof.} The formulas for $U_{n,c}(0)$, $U_{n,c}(-n)$,
$U_{n,c}(1)$ and $U_{n,c}(-n-1)$ in Section~\ref{q-sums}, \eqref{s3}--\eqref{s6}, imply that
the denominators from Lemma~\ref{recursion} are sums of at most
$3$ products. A lengthy but straightforward calculation shows that
theses denominators are actually single products. \qed

\medskip

We finally give the formulas for $L_{n,c,p}$,
$M_{n,c,p}$ and $G_{n,c,p}$.
\begin{lem}
\label{final}
The generating function $G_{n,c,p}$ is equal to
\begin{equation*}
G_{n,c,p}= q^{\binom{p+1}{2}} \qbinom{n}{p}
\begin{cases}
\frac{1}{[c+p;q]_{n+1}} \prod\limits_{i=0}^{n}
\frac{[c+2i;q]_{n-i+1}}{[2+2i;q]_{n-i}} & 2|c \\
\prod\limits_{i=1}^{n} \frac{[c+2i-1;q]_{n-i+1}}{[2i;q]_{n-i+1}} &
2 \not| c
\end{cases}.
\end{equation*}
For $L_{n,c,p}$ we have
\begin{equation*}
L_{n,c,p} =
\begin{cases}
\prod\limits_{i=1}^{n-1} \frac{[c+2 i + 1;q]_{n-i}}{[2 i;q]_{n-i}
[2 i;q]} \frac{[c+1;q]_{n-1} [1;q]_{n-1}}{2} \left(
\frac{q^{\binom{p+1}{2}} {\qbinom{n-1}{p}} [c;q]}{[c+p;q]_n} -
\frac{q^{\binom{p}{2}} {\qbinom{n-1}{p-1}} [c+2n;q]}{[c+p+1;q]_n}
\right) & 2|c \\
\prod\limits_{i=1}^{n-1} \frac{[c+2 i;q]_{n-i}}{[2 i;q]_{n-i} [2
i;q]} \frac{[1;q]_{n-1}}{2} \left( q^{\binom{p+1}{2}}
{\qbinom{n-1}{p}} - q^{\binom{p}{2}} {\qbinom{n-1}{p-1}} \right)
& 2 \not| c
\end{cases}
\end{equation*}
and for $M_{n,c,p}$ we have
\begin{multline*}
M_{n,c,p}= \frac{(-1)^{n-1} q^{(n-1)(2c+n)/2}}{[1;q]_{n-2}} \\
\times
\begin{cases}
\prod\limits_{i=1}^{n-1} \frac{[c+2i;q]_{n-i}}{[2i;q]_{n-i}}
\left( \frac{q^{\binom{p+1}{2}} {\qbinom{n-1}{p}}
[c;q]}{[c+p;q]_n} \left( \frac{1}{[n-1;q]} -
\frac{[c+2n-1;q]}{[c+n;q] [2n-2;q]} \right)  +
\frac{q^{\binom{p}{2}} {\qbinom{n-1}{p-1}}}{[c+p+1;q]_n}
\frac{[c+2n-1;q]_2}{[c+n;q] [2n-2;q]} \right) & 2|c \\
\prod\limits_{i=1}^{n-1} \frac{[c+2i+1;q]_{n-i-1}}{[2i;q]_{n-i}}
\frac{1}{[2n-2;q]} \left( q^{\binom{p+1}{2}+n-1} {\qbinom{n-1}{p}}
+ q^{\binom{p}{2}} {\qbinom{n-1}{p-1}} \right) & 2 \not| c
\end{cases}.
\end{multline*}
\end{lem}

{\it Proof.} We show the assertion by induction with respect to
$n$. For $n=1$ observe that $L_{1,c,0}=1/2$, $L_{1,c,1}=-1/2$ and
$M_{1,c,0}=M_{1,c,1}=1/2$. Moreover
\begin{equation*}
G_{1,c,0}=
\begin{cases}
1 +q^2 + q^4 + \ldots + q^c = \frac{1-q^{2+c}}{1-q^2} & \text{if
$c$ is even} \\
1 + q^2 + q^4 + \ldots + q^{c-1} = \frac{1-q^{1+c}}{1-q^2} &
\text{if $c$ is odd}
\end{cases}
\end{equation*}
and
\begin{equation*}
G_{1,c,1}=
\begin{cases}
q +q^3 + q^5 + \ldots + q^{c-1} = \frac{q(1-q^{c})}{1-q^2} &
\text{if $c$ is even} \\
q + q^3 + q^5 + \ldots + q^{c} = \frac{q(1-q^{1+c})}{1-q^2} &
\text{if $c$ is odd}
\end{cases}.
\end{equation*}
Assume that the formulas are proved for $n-1$. Then the formula
for $L_{n,c,p}$ and $M_{n,c,p}$ can be checked by using the
recursions in Lemma~\ref{recursion}, \eqref{s3}--\eqref{s6} and the
formula for $G_{n-1,c,p}$ which is true by
the induction hypothesis. By \eqref{s2} and \eqref{s1} we have
\begin{multline*}
G_{n,c,p} = \sum_{k=0}^c L_{n,c,p} U_{n,c}(k) q^k + M_{n,c,p} W_{n,c} (k) q^k \\
= L_{n,c,p} \frac{\left( (1+(-1)^{c}) [(c+2)/2;q^{2}]_{n-1} (1+q^{c+n}) + (1-(-1)^{c})
  [(c+1)/2;q^{2}]_{n} (1-q^{2})\right)}{[1/2;q^{2}]_{n-1} (1+q^{n-1}) (1+q^{n})}
 \\ + M_{n,c,p}
(-1)^{n-1} q^{(-n+1)(2c+n)/2} \frac{[1;q]^2_{n-1}
[c+1;q]_{2n-1}}{[1;q]_{2n-1}}.
\end{multline*}
A lengthy but straightforward calculation proves the formula for $G_{n,c,p}$. \qed

\medskip

Now we are able to explain why Theorem~\ref{main} implies
Krattenthaler's and the author's refinement of the Bender-Knuth
(ex-)Conjecture. Krattenthaler's refinement, see
\cite[Theorem~21]{kratt}, is the generating function $G_{n,c,p}$
we have computed in Lemma~\ref{final} and thus we have reproved
his result with different methods.

The author's refinement, see \cite[Theorem~1]{fischer}, is the generating function of strict plane partitions
with parts in $\{1,2,\ldots,n\}$, at most $c$ columns and $k$ parts equal to
$n$, i.e. the sum over all $p$'s, $0 \le p \le n$, of the generating function
in Theorem~\ref{main}. In order to deduce Theorem~1 in \cite{fischer} from
Theorem~\ref{main} of the present paper one has to show that
$$
\sum_{p=0}^{n} L_{n,c,p} =0
$$
and
$$ \sum_{p=0}^{n} (-1)^{n-1} q^{(n-1)(k-c) - \binom{n}{2})+k} M_{n,c,p}=
\frac{q^{k n}}{[1;q]_{n-1}} \prod_{i=1}^{n-1}
  \frac{[c+i+1;q]_{i-1}}{[i;q]_{i}},
$$
where
$$[k-c-n+1;q]_{n-1}=(-1)^{n-1} q^{(n-1)(k-c)-\binom{n}{2}} [1+c-k;q]_{n-1}
$$
explains the factor in front of $M_{n,c,p}$. However, Theorem~1 from
\cite{fischer} was proved with methods similar to the methods we
have used to prove Theorem~\ref{main} and thus we omit to show this
implication, because it is surely a detour to prove Theorem~1 from
\cite{fischer} in this way.

\section{Some basic hypergeometric identities}
\label{q-sums}

In this section we derive some basic hypergeometric identities
which were needed above. The notation is adopted from
\cite[page 1--6]{gasper}. In particular the basic hypergeometric
series is defined by
\begin{equation*}
{} _{r} \phi _{s} \!\left[
                       \begin{matrix}
                          a_1,\dots, a_r  \\
                          b_1,\dots, b_s
                       \end{matrix}; {\displaystyle q,  z} \right]=
\sum_{n=0}^{\infty} \frac{(a_1;q)_{n} \cdots
(a_r;q)_{n}}{(q;q)_{n} (b_1;q)_n
  \cdots (b_s;q)_{n}} \left( (-1)^{n} q^{\binom{n}{2}} \right)^{s-r+1} z^{n},
\end{equation*}
where the rising $q$-factorial $(a;q)_{n}$ is given by
$(a;q)_{n}:=\prod_{i=0}^{n-1} (1-a q^{i})$. (Observe that
$[x;q]_{n}=(q^{x};q)_{n}/(1-q)^{n}$.) All identities in this
section were handled with Krattenthaler's Mathematica package HYPQ
\cite{hypq}. A Mathematica-file containing the computations can be
downloaded from my webpage (\url{http://www.uni-klu.ac.at/~ifischer/}).

The list of identities is the following. (For the definitions of
$U_{n,c}(k)$ and $W_{n,c}(k)$ see Section~\ref{compute}.)
\begin{equation}
\label{s2}
 \sum_{k=0}^{c} W_{n,c}(k) q^{k} =(-1)^{n-1} q^{(-n+1)(2c+n)/2}
\frac{[1;q]^2_{n-1} [c+1;q]_{2n-1}}{[1;q]_{2n-1}}
\end{equation}
\begin{equation}
\label{s3} U_{n,c}(0) = \frac{2  \prod_{i=1}^{n-2}
(1+q^{i})}{(1+q)^{n-1}}
\begin{cases}
\frac{[c+1;q]_{n-1}}{[(c+1)/2;q]_{n-1}} & \text{if $c$ is even} \\
\frac{[c;q]_{n}}{[c/2;q]_{n} (1+q)} & \text{if $c$ is odd}
\end{cases}
\end{equation}
\begin{equation}
\label{s4} U_{n,c}(-n) = \frac{2   (1+q)^{n-1} \prod_{i=1}^{n-2}
(1+q^{i})}{q^{(n-1)(3n-2)/2}}
\begin{cases}
\frac{[(c+2)/2;q]_{n-1}}{[c+2;q]_{n-1}} & \text{if $c$ is even} \\
\frac{[(c+1)/2;q]_{n} (1+q)}{[c+1;q]_{n}} & \text{if $c$ is odd}
\end{cases}
\end{equation}
\begin{multline}
\label{s5} U_{n,c}(1) = \frac{2 q^n  \prod_{i=1}^{n-2}
(1+q^i)}{(1+q)^{n-1}} \\ \times
\begin{cases}
\frac{[c;q]_{n-1} [n-2;q]}{[(c+1)/2;q^2]_{n-1}}  & \text{if $c$ is even} \\
\frac{[c-1;q]_{n-1} [n-2;q]}{[c/2,q^2]_{n-1}} - q^{c+n-3}
\frac{[c;q]_{n-2} [2 \lfloor (n-1)/2 \rfloor +1;q] [2 \lfloor n/2
\rfloor;q]}{[(c+2)/2;q^2]_{n-1}}  & \text{if $c$ is odd}
\end{cases}
\end{multline}
\begin{multline}
\label{s6} U_{n,c}(-n-1) = \frac{2 (-1)^n (1+q)^{n-1}
\prod_{i=1}^{n-2} (1+q^i)}{q^{(n+1)(3n-4)/2}} \\ \times
\begin{cases}
\frac{[(c+4)/2;q^2]_{n-1} [n-2;q]}{[c+3;q]_{n-1}}  & \text{if $c$ is even} \\
\frac{[(c+3)/2;q^2]_{n-1} [n-2;q]}{[c+2;q]_{n-1}} + q^{c+n}
\frac{[(c+3)/2;q^2]_{n-1} [2 \lfloor (n-1)/2 \rfloor +1;q] [2
\lfloor n/2 \rfloor;q]}{[c+2;q]_{n} (1+q)^2}  & \text{if $c$ is
odd}
\end{cases}
\end{multline}
\begin{multline}
\label{s1}
\sum_{k=0}^{c} U_{n,c}(k) q^{k}  \\
 =  \frac{2}{[1/2;q^{2}]_{n-1} (1+q^{n-1}) (1+q^{n})}
\begin{cases}
[(c+2)/2;q^{2}]_{n-1} (1+q^{c+n}) & \text{if $c$ is even} \\
[(c+1)/2;q^{2}]_{n} (1-q^{2}) & \text{if $c$ is odd}
\end{cases}
\end{multline}

We first consider \eqref{s2}. Using the basic hypergeometric
notation introduced above, the left-hand-side can be written as
\begin{equation*}
{} _{2} \phi _{1} \! \left [                \begin{matrix} \let
\over / \def\frac\#1\#2{\#1 / \#2}
   q^ {n}, q^{-c} \\ \let \over / \def\frac\#1\#2{\#1 / \#2} q^ {1 - c -
   n}\end{matrix} ;q, {\displaystyle q} \right ]\,
 \frac{ ({\let \over / \def\frac\#1\#2{\#1 / \#2} q}; q) _{-1 + n} \,
  ({\let \over / \def\frac\#1\#2{\#1 / \#2} q^ {1 - c - n}}; q) _{-1 + n}}{(1-q)^{2n-2}}.
\end{equation*}
If we use the $q$-Chu-Vandermonde summation formula \cite[(1.5.3);
Appendix (II.6)]{gasper}
\begin{equation}
\label{S2101} {} _{2} \phi _{1} \! \left [
\begin{matrix} \let \over / \def\frac\#1\#2{\#1 / \#2}
   a, q^ {-n}\\ \let \over / \def\frac\#1\#2{\#1 / \#2} c \end{matrix} ;q,
 {\displaystyle q} \right ]\, = \frac{a^{n} (c/a;q)_{n}}{(c;q)_{n}},
\end{equation}
we obtain
\begin{equation}
\frac{q^{c\,n}\,({\let \over / \def\frac\#1\#2{\#1 / \#2} q}; q)
_{-1 + n} \,
    ({\let \over / \def\frac\#1\#2{\#1 / \#2} q^ {1 - c - 2\,n}}; q) _{c} \,
    ({\let \over / \def\frac\#1\#2{\#1 / \#2} q^ {1 - c - n}}; q) _{-1 + n}
  }{(1-q)^{2n-2} ({\let \over /
    \def\frac\#1\#2{\#1 / \#2} q^ {1 - c - n}}; q) _{c} }
\end{equation}
and this is equivalent to the right-hand-side in \eqref{s2}.

Next we consider \eqref{s3}, \eqref{s4}, \eqref{s5} and
\eqref{s6}. In all these summations we first use some contiguous
relations before we apply the following summation formula
\begin{equation}
\label{S2104} {} _{2} \phi _{1} \! \left [
\begin{matrix} \let \over / \def\frac\#1\#2{\#1 / \#2}
   a, b \\ \let \over / \def\frac\#1\#2{\#1 / \#2} a q/b \end{matrix} ;q,
 {\displaystyle -q/b} \right ]\, = \frac{(-q;q)_{\infty} (a q;q^{2})_{\infty}
 (a q^{2}/b^{2};q^{2})_{\infty}}{(-q/b;q)_{\infty} (a q/b;q)_{\infty}},
\end{equation}
see \cite[(1.8.1); Appendix (II.9)]{gasper}, where
$(a;q)_{\infty}=\prod_{i=0}^{\infty}(1-a q^{i})$. We define
$$
X_{n,c}(k)  = \sum_{i=1}^{n-1}   q^{\binom{i}{2}}
\frac{[k+1;q]_{n-1}
  [k-c-i+1;q]_{i-1} [k-c-n+1;q]_{n-i-1}}{[1;q]_{i-1} [1;q]_{n-1-i}
  [c+i+1;q]_{n-1}}
$$
and
$$
Y_{n,c}(k)  = \sum_{i=1}^{n-1} q^{\binom{i}{2}}
\frac{[k+1;q]_{i-1}
  [k+i+1;q]_{n-i-1} [k-c-n+1;q]_{n-1}}{ [1;q]_{i-1}
  [1;q]_{n-1-i} [c+i+1;q]_{n-1}}.
$$
Observe that
$$
U_{n,c}(k) = (-1)^{n} q^{(n-1)(2c+n)/2} ((-1)^{c} X_{n,c}(k) - Y_{n,c}(k)) +
(-1)^{k} q^{(n-1)k}.
$$
In order to prove \eqref{s3}, it suffices to compute $X_{n,c}(0)$,
$Y_{n,c}(0)$. Using basic hypergeometric notation $X_{n,c}(0)$ is
equal to
\begin{equation*}
{} _{2} \phi _{1} \! \left [                \begin{matrix} \let
\over / \def\frac\#1\#2{\#1 / \#2}
   q^{c+1}, q^{-n+2} \\ \let \over / \def\frac\#1\#2{\#1 / \#2} q^{c+n+1} \end{matrix} ;q,
 {\displaystyle -q^{n}} \right ]\, \frac{(q;q)_{n-1}
 (q^{1-c-n};q)_{n-2}}{(q;q)_{n-2} (q^{c+2};q)_{n-1}}.
\end{equation*}
If we apply the following contiguous relation
\begin{equation*}
{} _{r} \phi _{s} \! \left [                \begin{matrix} \let
\over / \def\frac\#1\#2{\#1 / \#2}
   a, (A) \\ \let \over / \def\frac\#1\#2{\#1 / \#2} b, (B) \end{matrix} ;q,
 {\displaystyle z} \right ]  =
\frac{1-b/q}{a-b/q} {} _{r} \phi _{s} \! \left [
\begin{matrix} \let \over / \def\frac\#1\#2{\#1 / \#2}
   a, (A) \\ \let \over / \def\frac\#1\#2{\#1 / \#2} b/q, (B) \end{matrix} ;q,
 {\displaystyle z/q} \right ]  +
\frac{1-a}{b/q-a} {} _{r} \phi _{s} \! \left [
\begin{matrix} \let \over / \def\frac\#1\#2{\#1 / \#2}
   a q, (A) \\ \let \over / \def\frac\#1\#2{\#1 / \#2} b, (B) \end{matrix} ;q,
 {\displaystyle z/q} \right ]
\end{equation*}
and then \eqref{S2104} we obtain the following  formula for
$X_{n,c}(0)$.
\begin{multline*}
\frac{q^{-1-c} (q^{1-c-n};q)_{n-2} (-q;q)_{\infty}
(q^{c+2};q^{2})_{\infty}
  (q^{c+2n-1};q^{2})_{\infty}}
{(q^{c+2};q)_{n-2} (-q^{n-1};q)_{\infty} (q^{c+n};q)_{\infty}} + \\
\frac{q^{-c-n} (q^{c+1};q)_{1} (q^{1-c-n};q)_{n-2} (q^{n-1};q)_{1}
  (-q;q)_{\infty} (q^{c+3};q^{2})_{\infty} (q^{c+2n};q^{2})_{\infty}}
{(q^{c+2};q)_{n-1} (q^{1-n};q)_{1} (-q^{n-1};q)_{\infty}
  (q^{c+n+1};q)_{\infty}}
\end{multline*}
Next observe that $Y_{n,c}(0)$ is equal to
\begin{equation*}
(-1)^{n-1}  {}_{3} \phi _{2} \! \left [
\begin{matrix} \let \over / \def\frac\#1\#2{\#1 / \#2}
   q^{-n+2}, q^{c+2}, q \\ \let \over / \def\frac\#1\#2{\#1 / \#2} q^{c+n+1}, q^{2} \end{matrix} ;q,
 {\displaystyle -q^{n-1}} \right ]\, \frac{(q^{n-1};q)_{1}
 (q^{c+1};q)_{1}}{(q;q)_{1} (q^{c+n};q)_{1}}.
\end{equation*}
Using the contiguous relation
\begin{equation}
\label{C02} {} _{r} \phi _{s} \! \left [
\begin{matrix} \let \over / \def\frac\#1\#2{\#1 / \#2}
   (A),q \\ \let \over / \def\frac\#1\#2{\#1 / \#2} (B) \end{matrix} ;q,
 {\displaystyle z} \right ] =
\frac{(-1)^{r+s} q^{1-r+s}}{z} \frac{\prod_{i=1}^{s}  (1-
  B_{i}/q)}{\prod_{i=1}^{r} (1- A_{i}/q)} \left(1 - {} _{r} \phi _{s} \! \left [
  \begin{matrix} \let \over / \def\frac\#1\#2{\#1 / \#2}
   (A/q),q \\ \let \over / \def\frac\#1\#2{\#1 / \#2} (B/q) \end{matrix} ;q,
 {\displaystyle q^{-1+r-s} z} \right ] \right),
\end{equation}
we transform the ${}_{3} \phi_{2}$ series into a ${}_{2} \phi_{1}$
series. Next we apply the following contiguous relation
\begin{equation}
\label{C15} {} _{r} \phi _{s} \! \left [
\begin{matrix} \let \over / \def\frac\#1\#2{\#1 / \#2}
   a, (A) \\ \let \over / \def\frac\#1\#2{\#1 / \#2} (B) \end{matrix} ;q,
 {\displaystyle z} \right ] =
{} _{r} \phi _{s} \! \left [                \begin{matrix} \let
\over / \def\frac\#1\#2{\#1 / \#2}
   a q, (A) \\ \let \over / \def\frac\#1\#2{\#1 / \#2} (B) \end{matrix} ;q,
 {\displaystyle z} \right ] + (-1)^{r+s} a z \frac{\prod_{i=1}^{r-1}
 (1-A_{i})}{\prod_{i=1}^{s} (1-B_{i})}
{} _{r} \phi _{s} \! \left [                \begin{matrix} \let
\over / \def\frac\#1\#2{\#1 / \#2}
   a q, (q A) \\ \let \over / \def\frac\#1\#2{\#1 / \#2} (q B) \end{matrix} ;q,
 {\displaystyle q^{1-r+s} z} \right ]
\end{equation}
and finally we are able to apply \eqref{S2104} to the remaining
two ${}_{2} \phi_{1}$ series. We obtain a formula for $Y_{n,c}(0)$
in terms of rising $q$-factorials. If we  combine this with the
formula for $X_{n,c}(0)$ we obtain \eqref{s3}.

The situation is similar for \eqref{s4}, \eqref{s5}, \eqref{s6},
however, we do not give the proofs in detail. In order to describe
a set of contiguous relations which is needed before \eqref{S2104}
can be applied in each of these cases, we use Krattenthaler's
table of contiguous relations, which can be found in the HYPQ
documentation \cite{hypq}. For details see the Mathematica-file
with the computation on my webpage
\url{http://www.uni-klu.ac.at/~ifischer/}.

\bigskip

\begin{tabular}{|l||l|} \hline
$X_{n,c}(-n)$  & C34 \\ \hline $Y_{n,c}(-n)$  & C15 \\ \hline
$X_{n,c}(1)$   & C14, C42, C34 \\ \hline $Y_{n,c}(1)$   & C16,
C15, C02, C41, C11, C37 \\ \hline $X_{n,c}(-n-1)$ & C14, C42, C34
\\ \hline $Y_{n,c}(-n-1)$ & C14, C15, C36, C11, C37 \\ \hline
\end{tabular}

\bigskip

Finally we consider \eqref{s1}, which is a double sum and thus
the most-complicated identity. First we compute $\sum_{k=0}^{c}
X_{n,c}(k) q^{k}$. The trick is to consider a more general expression.
Observe that $X_{n,c}(k)$ is the unique
$q$-polynomial in $k$ of degree at most $2n-3$ with the property
that $X_{n,c}(k)=0$ for $k=-1,-2,\ldots,-n+1$ and
$X_{n,c}(c+i)=(-1)^{n-1} q^{-\binom{n}{2}} (-q^{n-1})^{i}$ for
$i=1,2,\ldots,n-1$. (In Lemma~\ref{lagrange} $X_{n,c}(k)$ is
actually constructed in such a way that these conditions are
fulfilled.) Consequently $S_{n,c}(d): = \sum_{k=0}^{d} X_{n,c}(k)
q^{k}$ is the unique $q$-polynomial in $d$ of degree at most
$2n-2$ in $d$ with $S_{n,c}(d)=0$ for $d=-1,-2,\ldots,-n$ and
$$ S_{n,c}(c+i)=S_{n,c}(c)+ (-1)^{n-1} q^{-\binom{n}{2}+c+n}
\frac{-1+(-q^{n})^{i}}{1+q^{n}} $$ for $i=1,\ldots,n-1$. Thus, by
$q$-Lagrange interpolation,
\begin{multline*}
S_{n,c}(d) = \sum_{i=1}^{n} \left( S_{n,c}(c) + (-1)^{n-1}
q^{-\binom{n}{2}+c+n} \frac{(-q^{n})^{i} -
  1}{1+q^{n}} \right) q^{\binom{n-i+1}{2}} (-1)^{n+i} \\
\times \frac{[d+1;q]_{n} [d-c-i+2;q]_{i-1}
[d-c-n+1;q]_{n-i}}{[c+i;q]_{n} [1;q]_{i-1} [1;q]_{n-i}}.
\end{multline*}
Apriori the degree of this $q$-polynomial is $2n-1$. Thus the
coefficient of $(q^{d})^{2n-1}$
\begin{multline*}
 \sum_{i=1}^{n}
\left( S_{n,c}(c) + (-1)^{n-1} q^{-\binom{n}{2}+c+n}
\frac{(-q^{n})^{i} -
  1}{1+q^{n}} \right)
\frac{ q^{\binom{n-i+1}{2}} (-1)^{n+i} q^{-1+c+i+n-c
n}}{[c+i;q]_{n} [1;q]_{i-1} [1;q]_{n-i}},
\end{multline*}
must be zero. We obtain the following expression for
$S_{n,c}(c)=\sum_{k=0}^{c} X_{n,c}(k) q^{k}$.
$$
\frac{(-1)^{n-1} q^{-\binom{n}{2}+c+n}}{1+q^{n}} \left(
q^{\binom{n}{2}} \frac{\sum\limits_{i=1}^{n}
    \frac{q^{\binom{i+1}{2}}}{[c+i;q]_{n} [1;q]_{i-1} [1;q]_{n-i}}}
        {\sum\limits_{i=1}^{n} \frac{(-1)^{i} q^{\binom{n-i+1}{2}+i}}
                             {[c+i;q]_{n} [1;q]_{i-1} [1;q]_{n-i}}}+1 \right)
$$
This formula simplifies since
\begin{equation}
\label{hypo}
\sum_{i=1}^{n} \frac{(-1)^{i} q^{\binom{n-i+1}{2}+i}}
                             {[c+i;q]_{n} [1;q]_{i-1} [1;q]_{n-i}} =
-q^{(n-1) c + (n+1) n/2} \frac{[1;q]_{2n-2}}{[1;q]^{2}_{n-1}
[c+1;q]_{2n-1}}.
\end{equation}
In order to see that observe that the left-hand-side in this
equation is equal to
\begin{equation*}
-{} _{2} \phi _{1} \! \left [                \begin{matrix} \let
\over / \def\frac\#1\#2{\#1 / \#2}
   q^{c+1}, q^{-n+1} \\ \let \over / \def\frac\#1\#2{\#1 / \#2} q^{c+n+1} \end{matrix} ;q,
 {\displaystyle q} \right ]\, \frac{q^{n^{2}/2-n/2+1} (1-q)^{2n-1}}{(q;q)_{n-1} (q^{c+1};q)_{n}}.
\end{equation*}
Using \eqref{S2101} we obtain \eqref{hypo}. Thus
\begin{multline*}
S_{n,c}(c) = \frac{(-1)^{n-1} q^{-\binom{n}{2} + c + n}}{1+q^{n}}
\\ \times  \left( - q^{c- n - c n} \sum_{i=0}^{n-1} q^{\binom{i+2}{2}}
\frac{[c+1;q]_{i} [c +i+n+1;q]_{n-i-1} [n-i;q]_{i}}{[n;q]_{n-1}
[1;q]_{i}} + 1 \right).
\end{multline*}
Similarly one can show that
\begin{multline*}
T_{n,c}(d) = \sum_{i=0}^{n-1} \bigg( T_{n,c}(0) + \frac{(-1)^{n}
q^{(n-1)(2c+n)/2} (-1+(-q^{-n})^{i})}{1+q^{n}} (-1)^{i+n}
q^{(i+i^{2}+ n + 2 c n + 2 i n + n^{2})/2} \\
\times \frac{[d-c-n;q]_{n} [d;q]_{i}
[d+i+1;q]_{n-i-1}}{[c+i+1;q]_{n} [1;q]_{i}
  [1;q]_{n-1-i}} \bigg),
\end{multline*}
where $T_{n,c}(d) := \sum_{k=d}^{c} Y_{n,c}(k) q^{k}$. Again we
have that $T_{n,c}(d)$ is a $q$-polynomial of degree at most
$2n-2$, however, the left-hand-side of the equation above is
apriori a $q$-polynomial of degree $2n-1$. Thus we obtain the
following formula for $T_{n,c}(0)=\sum_{k=0}^{c} Y_{n,c}(k) q^{k}$.
$$
\frac{(-1)^{n-1} q^{-(n-1)(2c+n)/2}}{(1+q^{n})} \left(
\frac{q^{\binom{n}{2}} \sum_{i=0}^{n-1}
\frac{q^{\binom{i}{2}}}{[c+i+1;q]_{n}
    [1;q]_{i} [1;q]_{n-i-1}}}
{\sum_{i=0}^{n-1} \frac{(-1)^{i} q^{\binom{n+i}{2}}}{[c+i+1;q]_{n}
[1;q]_{i}
    [1;q]_{n-i-1}}} - 1 \right)
$$
Again the formula simplifies since
\begin{equation}
\sum_{i=0}^{n-1} \frac{(-1)^{i} q^{\binom{n+i}{2}}}{[c+i+1;q]_{n}
[1;q]_{i}
  [1;q]_{n-i-1}} = \frac{q^{\binom{n}{2}} [1;q]_{2n-2}}{[c+1;q]_{2n-2}
  [1;q]^{2}_{n-1}}.
\end{equation}
In order to see that transform the sum into hypergeometric notation
\begin{equation*}
{} _{2} \phi _{1} \! \left [                \begin{matrix} \let
\over / \def\frac\#1\#2{\#1 / \#2}
   q^{c+1}, q^{1-n} \\ \let \over / \def\frac\#1\#2{\#1 / \#2} q^{c+n+1} \end{matrix} ;q,
 {\displaystyle q^{2n-1}} \right ]\, \frac{q^{n^{2}/2-n/2} (1-q)^{2n-1}}{(q;q)_{n-1} (q^{c+1};q)_{n}}
\end{equation*}
and apply the summation formula, see \cite[(1.5.2); Appendix
(II.7)]{gasper},
\begin{equation*}
\label{S2102} {} _{2} \phi _{1} \! \left [
\begin{matrix} \let \over / \def\frac\#1\#2{\#1 / \#2}
   a, q^{-n} \\ \let \over / \def\frac\#1\#2{\#1 / \#2} c \end{matrix} ;q,
 {\displaystyle \frac{c q^{n}}{a} } \right ]\, = \frac{(c/a;q)_{n}}{(c;q)_{n}}
\end{equation*}
to obtain the result. This implies that
$$
T_{n,c}(0) = \frac{(-1)^{n-1} q^{-(n-1)(2c+n)/2}}{(1+q^{n})}
\left( \sum_{i=0}^{n-1} q^{\binom{i}{2}} \frac{[c+1;q]_{i}
[c+i+n+1;q]_{n-i-1}
  [n-i;q]_{i}}{[n;q]_{n-1} [1;q]_{i}} - 1 \right).
$$
Therefore
\begin{multline*}
\sum_{k=0}^c U_{n,c}(k) q^k =\frac{1+(-1)^{c} q^{(1+c) n}}{1+q^{n}} +  \frac{(-1)^n
q^{-(n-1)(2c+n)/2}}{1+q^n} \\ \times \bigg( \sum_{i=0}^{n-1} (-1)^c
q^{\binom{i+2}{2}+c} \frac{[c+1;q]_i [c+i+n+1;q]_{n-i-1}
[n-i;q]_i}{[n;q]_{n-1} [1;q]_i} \\ -
\sum_{i=0}^{n-1} q^{\binom{i}{2}} \frac{[c+1;q]_i
[c+i+n+1;q]_{n-i-1} [n-i;q]_i}{[n;q]_{n-1} [1;q]_i} - (-q^n)^{c+1}
+ 1 \bigg)
\end{multline*}
Consequently it suffices to compute
\begin{multline}
\label{e1} \sum_{i=0}^{n-1} q^{\binom{i+2}{2}+c} \frac{[c+1;q]_i
[c+i+n+1;q]_{n-i-1} [n-i;q]_i}{[n;q]_{n-1} [1;q]_i} \\ -
\sum_{i=0}^{n-1} q^{\binom{i}{2}} \frac{[c+1;q]_i
[c+i+n+1;q]_{n-i-1} [n-i;q]_i}{[n;q]_{n-1} [1;q]_i}
\end{multline}
and
\begin{multline}
\label{e2} \sum_{i=0}^{n-1} q^{\binom{i+2}{2}+c} \frac{[c+1;q]_i
[c+i+n+1;q]_{n-i-1} [n-i;q]_i}{[n;q]_{n-1} [1;q]_i} \\ +
\sum_{i=0}^{n-1} q^{\binom{i}{2}} \frac{[c+1;q]_i
[c+i+n+1;q]_{n-i-1} [n-i;q]_i}{[n;q]_{n-1} [1;q]_i}.
\end{multline}
Using basic hypergeometric notation \eqref{e1} is equal to
\begin{equation*}
 {} _{4} \phi _{3} \! \left [
\begin{matrix} \let \over / \def\frac\#1\#2{\#1 / \#2}
   q^{1+c}, q^{3/2+c/2},-q^{3/2+c/2},q^{1-n} \\ \let \over / \def\frac\#1\#2{\#1 / \#2}
   q^{1/2+c/2}, -q^{1/2+c/2}, q^{1+c+n} \end{matrix} ;q,
 {\displaystyle -q^{n-1}} \right ]\, \frac{(q^{1+c};q)_1
 (q^{1+c+n};q)_{n-1}}{(q^n;q)_{n-1}}.
\end{equation*}
We apply the following transformation
\begin{equation}
\label{T4308}
 {} _{4} \phi _{3} \! \left [
\begin{matrix} \let \over / \def\frac\#1\#2{\#1 / \#2}
   a, a^{1/2} q ,- a^{1/2} q , b \\ \let \over / \def\frac\#1\#2{\#1 / \#2}
   a^{1/2}, -a^{1/2}, a q/b \end{matrix} ;q,
 {\displaystyle t} \right ]\, =
{} _{2} \phi _{1} \! \left [
\begin{matrix} \let \over / \def\frac\#1\#2{\#1 / \#2}
   1/b, t \\ \let \over / \def\frac\#1\#2{\#1 / \#2}
   b q t \end{matrix} ;q,
 {\displaystyle a q} \right ]\, \frac{(a q;q)_{\infty} (b
 t;q)_{\infty}}{(t;q)_{\infty} (a q/b;q)_{\infty} },
\end{equation}
which can be found in \cite[Ex. 2.2]{gasper} and obtain a ${}_2
\phi_{1}$-series.  We apply another transformation
\begin{equation}
\label{T2101} {} _{2} \phi _{1} \! \left [
\begin{matrix} \let \over / \def\frac\#1\#2{\#1 / \#2}
   a, b \\ \let \over / \def\frac\#1\#2{\#1 / \#2}
   c \end{matrix} ;q,
 {\displaystyle z} \right ]  =
 {} _{2} \phi _{1} \! \left [
\begin{matrix} \let \over / \def\frac\#1\#2{\#1 / \#2}
   c/b, z \\ \let \over / \def\frac\#1\#2{\#1 / \#2}
   a z \end{matrix} ;q,
 {\displaystyle b} \right ] \frac{(b;q)_{\infty} (a z;q)_{\infty}}{(c;q)_{\infty}
(z;q)_{\infty}} ,
\end{equation}
see \cite[(1.4.1); Appendix (III.1)]{gasper}, before we are able
to apply the summation \eqref{S2104}. In basic hypergeometric
notation \eqref{e2} is equal to
\begin{equation*}
 {} _{4} \phi _{3} \! \left [
\begin{matrix} \let \over / \def\frac\#1\#2{\#1 / \#2}
   q^{1+c}, i q^{3/2+c/2},- i q^{3/2+c/2},q^{1-n} \\ \let \over / \def\frac\#1\#2{\#1 / \#2}
   i q^{1/2+c/2}, -i q^{1/2+c/2}, q^{1+c+n} \end{matrix} ;q,
 {\displaystyle -q^{n-1}} \right ]\, \frac{(q^{2+2c};q)_1
 (q^{1+c+n};q)_{n-1}}{(q^{c+1};q)_{1} (q^n;q)_{n-1}}.
\end{equation*}
We apply the following transformation rule {\small
\begin{multline}
\label{T4310}
 {} _{4} \phi _{3} \! \left [
\begin{matrix} \let \over / \def\frac\#1\#2{\#1 / \#2}
   a, b, c, d \\ \let \over / \def\frac\#1\#2{\#1 / \#2}
   a q/b, a q/c, a q/d \end{matrix} ;q,
{\displaystyle -a q^{2}/(bcd)} \right ]\, =
\\
 {} _{8} \phi _{7} \! \left [
\begin{matrix} \let \over / \def\frac\#1\#2{\#1 / \#2}
   a^{2} q/(b c d), a q^{3/2}/(b c d)^{1/2},  -a q^{3/2}/(b c d)^{1/2},
   a^{1/2},-a^{1/2}, aq /(cd), aq / (bd), aq /(bc)  \\ \let \over / \def\frac\#1\#2{\#1 / \#2}
   a q^{1/2}/(bcd)^{1/2}, - a q^{1/2} /(bcd)^{1/2}, a^{3/2} q^{2}/ (b c d),
   -a^{3/2} q^{2}/ (b c d),
aq/b, aq/c, aq/d  \end{matrix} ;q,
 {\displaystyle -q} \right ]\ \\
\times \frac{(aq;q)_{\infty} (-q;q)_{\infty} (a^{3/2} q^{2}/(b c
d);q)_{\infty}
 (-a^{3/2} q^{2}/(b c d);q)_{\infty}}
  {(a^{2} q^{2}/(b c d);q)_{\infty} (-a
  q^{2}/(b c d);q)_{\infty} (a^{1/2} q;q)_{\infty} (-a^{1/2}
  q;q)_{\infty}}
\end{multline}}
see \cite[Ex. 2.13 (ii)]{gasper}, to obtain a ${}_{8}
\phi_{7}$-series. Next we apply the transformation rule {\small
\begin{multline}
\label{T8704}
 {} _{8} \phi _{7} \! \left [
\begin{matrix} \let \over / \def\frac\#1\#2{\#1 / \#2}
   a, a^{1/2} q, -a^{1/2} q, b, c, d, e, f \\ \let \over / \def\frac\#1\#2{\#1 / \#2}
   a^{1/2}, -a^{1/2}, aq/b, aq/c, aq/d, aq/e, aq/f \end{matrix} ;q,
 {\displaystyle a^{2} q^{2}/(b c d e f)} \right ]\, = \\
 {} _{8} \phi _{7} \! \left [
\begin{matrix} \let \over / \def\frac\#1\#2{\#1 / \#2}
   a^{2} q/(b c d), a q^{3/2}/(b c d)^{1/2},  -a q^{3/2}/(b c d)^{1/2},
   a q/(c d), a q/(b d), a q/(b c), e, f  \\ \let \over / \def\frac\#1\#2{\#1 / \#2}
   a q^{1/2}/(b c d)^{1/2}, -a q^{1/2}/(b c d)^{1/2}, a q/b, a q/c, a q/d,
   a^{2} q^{2}/(b c d e), a^{2} q^{2}/(b c d f)  \end{matrix} ;q,
 {\displaystyle a q/(e f) } \right ]\, \\
 \times
 \frac{(aq;q)_{\infty} (aq/(e f);q)_{\infty} (a^{2} q^{2}/(b c d
e);q)_{\infty}
  (a^{2} q^{2}/(b c d f);q)_{\infty}}{(a q/e;q)_{\infty} (a q /f;q)_{\infty}
  (a^2 q^{2}/(b c d);q)_{\infty} (a^{2} q^{2}/(b c d e f);q)_{\infty}},
\end{multline}}
see \cite[(2.10.1); Appendix (III.23)]{gasper} and finally the
summation formula {\small \begin{multline} \label{S8702}
 {} _{8} \phi _{7} \! \left [
\begin{matrix} \let \over / \def\frac\#1\#2{\#1 / \#2}
   -(a b/q)^{1/2}c, i (a b)^{1/4} c^{1/2} q^{3/4},
   -i (a b)^{1/4} c^{1/2} q^{3/4}, a, b, c, -c, -(a b q)^{1/2}/c \\ \let \over / \def\frac\#1\#2{\#1 / \#2}
   i (a b/q)^{1/4} c^{1/2}, -i (a b/q)^{1/4} c^{1/2}, -
   (b q/a)^{1/2} c,-(a q/b)^{1/2} c,  -(a b q)^{1/2}, (a b q)^{1/2}, c^{2} \end{matrix} ;q,
 {\displaystyle \frac{c q^{1/2}}{(a b)^{1/2}}} \right ]\, = \\
\frac{(-(a b q)^{1/2} c;q)_{\infty} (-c (q/ab)^{1/2};q)_{\infty}}
{(-(b q/a)^{1/2} c;q)_{\infty} (-(a q/b)^{1/2} c;q)_{\infty}}
\frac{(aq;q^{2})_{\infty} (b q;q^{2})_{\infty} (c^{2}
q/a;q^{2})_{\infty}
  (c^{2} q/b;q^{2})_{\infty}}{(q;q^2)_{\infty} (a b q;q^2)_{\infty} (c^{2}
 q;q^2)_{\infty} (c^{2} q/(a b);q^2)_{\infty} },
\end{multline}}
see \cite[Ex. 2.17(i); Appendix (II.16)]{gasper}, in order to
obtain the closed form for \eqref{e2}.

\end{document}